\documentclass[12pt]{article}

\usepackage{amsmath,amssymb,amsfonts,amsthm}
\usepackage{mathrsfs,bm}
\usepackage{xcolor}
\usepackage[numbers,sort&compress]{natbib}
\usepackage{enumitem}
\usepackage{hyperref}
\usepackage{microtype}

\allowdisplaybreaks

\hypersetup{
    colorlinks=true,
    linkcolor=blue,
    citecolor=red,
    urlcolor=blue
}

\newtheorem{theorem}{Theorem}[section]
\newtheorem{lemma}[theorem]{Lemma}
\newtheorem{remark}[theorem]{Remark}
\newtheorem{proposition}[theorem]{Proposition}
\newtheorem{corollary}[theorem]{Corollary}
\numberwithin{equation}{section}

\newcommand{\R}{\mathbb{R}}
\newcommand{\eps}{\varepsilon}
\newcommand{\A}{\mathscr A_\alpha}
\newcommand{\HLS}{\mathcal{H}_\alpha}
\newcommand{\X}{X_\eps}
\newcommand{\M}{M_\eps}
\newcommand{\Xnorm}[1]{\|#1\|_{X_\eps}}
\DeclareMathOperator{\SL}{SL}
\newcommand{\Scale}[1]{\sigma_\eps(#1)}

\title{Sharp asymptotics and stability for a two-power affine Hardy-Littlewood-Sobolev inequality}
\author{Zhenzhen Wei, Denghui Wu$^{*}$}
\date{}

\begin{document}
\maketitle

\begingroup
\makeatletter
\renewcommand{\thefootnote}{}
\renewcommand{\@makefntext}[1]{\noindent#1}
\makeatother
\begin{NoHyper}
\footnotetext{


$^{*}$\textit{Corresponding author.}}
\end{NoHyper}
\endgroup

\begin{abstract}
We study the affine Hardy-Littlewood-Sobolev inequality under the constraint $\frac{1}{2}\int_{\mathbb{R}^n}(f^{q-\varepsilon}+f^{q+\varepsilon}) dx=1$ for $f\ge0$, where $q=\frac{2n}{n+\alpha}$, $n\ge2$, $0<\alpha<n$, and $0<\varepsilon<q-1$. We determine the second-order expansion of its sharp constant and the scale selected by all normalized radial extremals as $\varepsilon\rightarrow0^+$. For sufficiently small fixed $\varepsilon$, we prove stability estimates with power $2$ in $L^{q-\varepsilon}\cap L^{q+\varepsilon}$ for the Euclidean and affine energies. The power $2$ is sharp, and the sharp stability constants admit $O(\varepsilon^2)$ upper bounds as $\varepsilon\rightarrow0^+$. The proof combines identities for the two power integrals with Euclidean spectral estimates, affine compactness and a nonlinear directional expansion. No uniqueness or additional regularity of the radial extremals is assumed.
\end{abstract}

\section{Introduction}

The Hardy-Littlewood-Sobolev (HLS) inequality is a fundamental
estimate for Riesz potentials. Let $n\ge2$, $0<\alpha<n$, and
$q=\frac{2n}{n+\alpha}$. For $1\le r<\infty$, write
$\|f\|_r=(\int_{\mathbb{R}^n}|f|^r dx)^{\frac{1}{r}}$.
For real $f,h\in L^q(\mathbb{R}^n)$, the Riesz energy is
\[
\HLS(f,h)=\iint_{\R^n\times\R^n}
\frac{f(x)h(y)}{|x-y|^{n-\alpha}}dxdy.
\]
Lieb \cite{L83} obtained the sharp constant $C_0$ and characterized all equality
cases in
\[|\HLS(f,h)|\le C_0\|f\|_q\|h\|_q.\]
The classical Sobolev inequalities and their extensions in spherical and
logarithmic versions were studied by Talenti \cite{T76}, Beckner \cite{Bec93},
Frank and Lieb \cite{FL}, and Carlen and Loss \cite{CL}.
The fundamental tools in this theory are
rearrangement inequalities and their equality cases, see, e.g., \cite{BLL,Bur,BH,LL01}.

In the case $f=h$, the HLS inequality becomes
\[|\HLS(f,f)|\le C_0\|f\|_q^2.\]
This inequality is invariant under translations and dilations,
since translations $f_z(x)=f(x-z)$ and dilations $\mathcal D_\rho f(x)=\rho^{\frac{n}{q}}f(\rho x)$, $\rho>0$
preserve both $\HLS(f,f)$ and $\|f\|_q$.

To formulate the affine counterpart of the HLS inequality,
let us recall the affine Riesz energy derived from the directional correlations of $f$.
For $0\le f\in L^q$, the affine Riesz energy is defined by
\[
\A(f)=n\omega_n^{1-\frac{\alpha}{n}}
\bigg(\frac{1}{n}\int_{S^{n-1}}D_f(\xi)^{\frac{n}{\alpha}}d\xi\bigg)^{\frac{\alpha}{n}}.
\]
where
\[D_f(\xi)=\int_0^\infty t^{\alpha-1}\int_{\R^n}f(x)f(x+t\xi)dxdt,\quad \xi\in S^{n-1},\]
and $\omega_n$ is the volume of the Euclidean unit ball and
$d\xi$ is the surface area measure on the unit sphere $S^{n-1}=\{\xi\in\R^n:|\xi|=1\}$.
For a radial function $f(x)=F(|x|)$, $D_f(\xi)$ is independent of
$\xi$, then $\A(f)=\HLS(f,f)$. The energy $\A$ is invariant under
translations, the dilations $\mathcal D_\rho$, and the action of
$\SL(n)=\{T\in\R^{n\times n}:\det T=1\}$.

Haddad and Ludwig \cite{HL25b} established the affine refinement
\[
\HLS(f,f)\le\A(f)\le\HLS(f^\star,f^\star)\le C_0\|f\|_q^2.
\]
Here $f^\star$ is the symmetric decreasing rearrangement of $f$,
which is a nonnegative radial nonincreasing function
with satisfying $|\{f^\star>t\}|=|\{|f|>t\}|$ for $t>0$, and $|E|$ denotes
Lebesgue measure.
The classical affine Sobolev inequalities were established
by Zhang \cite{Zhang99} for $p=1$, Lutwak, Yang and Zhang \cite{LYZ} for $1<p<n$, Cianchi, Lutwak, Yang and Zhang \cite{CLYZ} for the limiting case $p = n$ and super-limiting case $p>n$.
The affine fractional Sobolev and isoperimetric inequalities were proved by Haddad and Ludwig \cite{HL25a}.
Hou, Xiao and Ye \cite{HXY18} studied the mixed volumes associated
with anisotropic Riesz potentials. Limits in the Riesz order were
investigated by Lin, Lan and Zhou \cite{LLZ26}.

Riesz energies under constraints involving two Lebesgue norms were
studied by Greco, Huang, Liu and Moroz \cite{GHLM}, who proved
existence of extremals when the two powers lie on opposite sides
of the HLS exponent. Weighted Riesz extremal problems on the
Heisenberg group were treated in \cite{CLT19}.
Interpolation constraints also occur in function-space inequalities
\cite{BS,DLLB22,DLLL22} and in weighted Sobolev and logarithmic
inequalities \cite{DDLL25,LL17,LL18}.
Related questions concerning sharp constants, extremal attainment
and norm constraints also arise in exponential inequalities
\cite{LL13,LLZ17,Lam17,LLZ19,LLZDN19}, including their affine
versions \cite{LLT,DLL}.

Quantitative stability asks how the deficit controls distance from
the extremal set. The Sobolev theorem of Bianchi and Egnell
\cite{BE} was extended to fractional inequalities in \cite{CFW},
and Carlen \cite{Car} developed a duality approach.
Explicit quantitative bounds were obtained in
\cite{DEFFL,CLT24,CLT25,CLT26}, with logarithmic stability studied
in \cite{CLT23} and stability for HLS critical points in
\cite{CLT26b}. Related weighted stability estimates appear in
\cite{CFLL,LLR26,DLL26}. Sharp local and global stability
constants need not coincide \cite{K23}.
Fan, Li and Zhang~\cite[Theorem 1.1]{FLZ} proved quadratic stability
for the affine fractional $L^2$-Sobolev inequality using a distance
defined by the $\dot H^s(\R^n)$ seminorm, $0<s<1$.
For the Riesz energies studied here, the two-power constraint
selects a dilation scale for the extremals, and stability is measured
in $L^{q-\eps}(\R^n)\cap L^{q+\eps}(\R^n)$.

Let $a=q-\eps$ and $b=q+\eps$, where
\begin{equation}\label{e1}
0<\eps<q-1.
\end{equation}
Thus $1<a<q<b$ and $a+b=2q$. Define
\begin{gather*}
\varphi_\eps(t)=\frac{1}{2}(t^a+t^b)\ \text{for}\ t\ge0,\quad
\M(f)=\int_{\R^n}\varphi_\eps(|f|)dx,\\*
\X=L^a(\R^n)\cap L^b(\R^n),\quad
\Xnorm{f}=\max\{\|f\|_a,\|f\|_b\}.
\end{gather*}
We call $f$ normalized if $\M(f)=1$. Since
$\varphi_\eps(t)\ge t^q$, normalization implies $\|f\|_q\le1$.
Interpolation also implies $\|h\|_q\le\Xnorm{h}$.
The two sharp constants are
\begin{align*}
C_{\eps,\rm E}&=\sup_{f\in\X,\ f\ge0,\ \M(f)=1}\HLS(f,f),\\
C_{\eps,\rm A}&=\sup_{f\in\X,\ f\ge0,\ \M(f)=1}\A(f).
\end{align*}
Rearrangement implies $C_{\eps,\rm E}=C_{\eps,\rm A}$.
Write their common value as $C_\eps$.
Minimizing the normalization factor along dilations yields
\[
C_\eps=\sup_{0\ne f\in\X,\ f\ge0}
\frac{\HLS(f,f)}{\big(\int_{\R^n}f^{q-\eps}dx
\int_{\R^n}f^{q+\eps}dx\big)^{\frac{1}{q}}}.
\]
The Euclidean existence result in
\cite[Theorem 1.2 and Remark 1.2]{GHLM} applies to this problem.
The role of the constraint can be seen from the dilation identities
\[
\int_{\R^n}|\mathcal D_\rho f|^{q-\eps}dx
=\rho^{-\frac{n\eps}{q}}\int_{\R^n}|f|^{q-\eps}dx,
\quad
\int_{\R^n}|\mathcal D_\rho f|^{q+\eps}dx
=\rho^{\frac{n\eps}{q}}\int_{\R^n}|f|^{q+\eps}dx.
\]
The energy is unchanged, while the two integrals scale in opposite
directions. As shown in Theorem~\ref{t4}, every normalized extremal
satisfies
\[
\int_{\R^n}u^a dx=\int_{\R^n}u^b dx=1.
\]
These equalities fix the dilation scale along each extremal orbit.
Let $\mathcal E_\eps$ be the nonempty compact set of normalized
radial nonincreasing Euclidean extremals in $X_\eps$, whose
compactness is also proved in Section~\ref{s2}.

These variational properties do not determine the loss in $C_\eps$
or the selected scale as $\eps\rightarrow0^+$. A further question is whether
the deficit controls distance in $X_\eps$ to the perturbed extremal
set. Classical $L^q$ stability concerns the classical HLS family,
which contains all dilations, and measures distance in a
single Lebesgue norm. The estimate sought here must incorporate the
perturbed extremals and control both powers in the constraint.
For the affine energy, it must also account for the freedom under
$\SL(n)$.

Our first result determines the second order loss and the
selected scale. Let $\Gamma(t)=\int_0^\infty s^{t-1}e^{-s}ds$,
$t>0$, be the Gamma function. Set
$\gamma=\frac{n+\alpha}{2}$ and
$Z_n=\frac{\pi^{\frac{n}{2}}\Gamma(\frac{n}{2})}{\Gamma(n)}$.
The normalized classical extremal is
\begin{equation}\label{e10}
U(x)=Z_n^{-\frac{1}{q}}(1+|x|^2)^{-\gamma},
\end{equation}
with $\|U\|_q=1$ and $\HLS(U,U)=C_0$.
Let $\psi(t)=\frac{\Gamma'(t)}{\Gamma(t)}$ and $\psi_1(t)=\psi'(t)$
be the digamma and trigamma functions. Define $\mathcal V_{n,\alpha}
=\gamma^2\big(\psi_1(\frac{n}{2})-\psi_1(n)\big),\
\kappa_{n,\alpha}=\frac{\mathcal V_{n,\alpha}}{q}>0$, $U_*=\mathcal D_{\rho_*}U$ and
\begin{align}
\rho_*&=Z_n^{\frac{1}{n}}
\exp\Big(\psi(n)-\psi\Big(\frac{n}{2}\Big)\Big).\label{e11}
\end{align}

\begin{theorem}\label{t1}
For fixed $n\ge2$ and $0<\alpha<n$,
\begin{equation}\label{e13}
C_\eps=C_0\big(1-\kappa_{n,\alpha}\eps^2+o(\eps^2)\big)
\ \text{as }\eps\rightarrow0^+.
\end{equation}
Every choice $u_\eps\in\mathcal E_\eps$ satisfies
\begin{gather}
\|u_\eps\|_q^2
=1-\kappa_{n,\alpha}\eps^2+o(\eps^2)
\ \text{as }\eps\rightarrow0^+,\label{e14}\\
u_\eps\rightarrow U_*\quad\text{strongly in }L^q(\R^n)
\ \text{as }\eps\rightarrow0^+.\label{e15}
\end{gather}
\end{theorem}

The second result concerns stability in the intersection norm.
Define the normalized extremal sets by
\begin{gather*}
\mathcal M_{\eps,{\rm E}}
=\{u(\cdot-y):u\in\mathcal E_\eps,\ y\in\R^n\},\\
\mathcal M_{\eps,{\rm A}}
=\{u(T(\cdot-y)):u\in\mathcal E_\eps,\ T\in\SL(n),\ y\in\R^n\}.
\end{gather*}
For normalized $f\ge0$, the Euclidean $X_\eps$ distance is
\[
d_{X_{\eps},\rm E}(f)
=\inf_{v\in\mathcal M_{\eps,{\rm E}}}\Xnorm{f-v}.
\]
The affine $X_\eps$ distance is
\[
d_{X_{\eps},{\rm A}}(f)
=\inf_{v\in\mathcal M_{\eps,{\rm A}}}\Xnorm{f-v}.
\]
The affine $L^q$ distance is
\[
d_{q,{\rm A}}(f)
=\inf_{v\in\mathcal M_{\eps,{\rm A}}}\|f-v\|_q.
\]
The Euclidean $L^q$ distance uses $\mathcal M_{\eps,{\rm E}}$
in the last infimum.

\begin{theorem}\label{t2}
There exists
$0<\eps_0<\min\{q-1,2-q,\frac{n}{\alpha}-q\}$,
depending on $n,\alpha$, such that, for every fixed
$0<\eps<\eps_0$, there are $k_{\eps,{\rm E}},k_{\eps,{\rm A}}>0$ with
\begin{align}
C_\eps-\HLS(f,f)
&\ge k_{\eps,{\rm E}}\big(d_{X_{\eps},\rm E}(f)\big)^2,\label{e18}\\
C_\eps-\A(f)&\ge k_{\eps,{\rm A}}d_{X_{\eps},{\rm A}}(f)^2\notag
\end{align}
for every $f\ge0$ with $\M(f)=1$.
Both estimates also hold with the corresponding $L^q$ distances.
The power $2$ is sharp in both estimates, whether the distance
is measured in $X_\eps$ or in $L^q$.
Moreover, $\mathcal E_\eps$ is finite.
\end{theorem}

The limiting scale in Theorem~\ref{t1} can also be characterized by
\[
\int_{\R^n}U_*^q\log(U_*)dx=0.
\]
For all $0<\eps<q-1$, Theorems~\ref{t8} and~\ref{t9} identify the
extremal sets and establish compactness of maximizing sequences
modulo translations and, in the affine case, volume-preserving
linear maps.
The stability proof does not assume uniqueness of the radial
extremal. Normalized dilations yield $O(\eps^2)$ upper bounds for
the sharp stability constants as $\eps\rightarrow0^+$.
A matching lower bound remains open.

The asymptotic proof combines classical HLS stability with the
identities in Theorem~\ref{t4} and a nonnegative logarithmic double integral. For Euclidean stability, two linearized constraints control the
constant and dilation directions in the spectral estimate.
For affine stability, Christ's quantitative rearrangement theorem
\cite{Chr} provides a geometric normalization, while directional
compactness and a nonlinear expansion control the affine
correction to second order.
Compactness uses concentration compactness \cite{L84,L85} and the
Br\'ezis-Lieb lemma \cite{BL}.
Section~\ref{s2} develops the variational properties,
Section~\ref{s3} proves the asymptotics, and Section~\ref{s4} treats
extremals and compactness. Sections~\ref{s5}-\ref{s6} establish
stability, and Section~\ref{s7} studies its sharp power and parameter
dependence.

\section{Preliminaries and variational properties}\label{s2}

Throughout, $0<\eps<q-1$, $a=q-\eps$ and $b=q+\eps$.
Unless another domain is specified, integrals are over $\R^n$.
For $R>0$ and $y\in\R^n$, let
$B_R(y)=\{x\in\R^n:|x-y|<R\}$, and write $B_R=B_R(0)$.
For $1<p<\infty$, $p'=\frac{p}{p-1}$ denotes the conjugate exponent.
For measurable sets $E,F$, their symmetric difference is
$E\triangle F=(E\setminus F)\cup(F\setminus E)$.
The notation $1_E$ denotes the indicator of $E$.
Constants may change from line to line, with dependence indicated
when needed.
The classical sharp constant is
\[
C_0=\pi^{\frac{n-\alpha}{2}}
\frac{\Gamma(\frac{\alpha}{2})}{\Gamma(\frac{n+\alpha}{2})}
\bigg(\frac{\Gamma(n)}{\Gamma(\frac{n}{2})}\bigg)^{\frac{\alpha}{n}}.
\]
For $f\ne0$, let $\Scale{f}$ be the unique $s>0$ such that
$\M(\frac{f}{s})=1$, and set $\Scale{0}=0$.

\begin{lemma}
For $0\ne f\in X_\eps$, the normalization factor $s=\Scale{f}$ solves
\begin{equation}\label{e20}
\frac{\int_{\R^n}|f|^a dx}{2s^a}+\frac{\int_{\R^n}|f|^b dx}{2s^b}=1
\end{equation}
and satisfies
\begin{gather}
\max\{2^{-\frac{1}{a}}\|f\|_a,2^{-\frac{1}{b}}\|f\|_b\}
\le\Scale{f}\le\Xnorm{f},\label{e21}\\
\|f\|_q\le\Scale{f}\le\Xnorm{f},\quad
\Scale{f^\star}=\Scale{f}.\notag
\end{gather}
\end{lemma}
\begin{proof}
The left side of \eqref{e20} decreases continuously from infinity
to zero. Its two terms yield the lower bound in \eqref{e21}.
Substitution of $s=\Xnorm{f}$ yields the upper bound.
Since $\varphi_\eps(t)\ge t^q$, normalization implies
$\|f\|_q\le\Scale{f}$. Equimeasurability proves the last assertion.
\end{proof}

A measurable star set $K$ is specified by its radial function
$\rho_K$. Its gauge satisfies $\|t\xi\|_K=\frac{t}{\rho_K(\xi)}$ and
$|K|=n^{-1}\int_{S^{n-1}}\rho_K^n d\xi$. Define $S_\alpha f$ by
$\rho_{S_\alpha f}^\alpha=D_f$.
Radial functions are understood up to surface-null sets, and star
sets up to Lebesgue-null sets. Positivity and finiteness of a radial
function below are required almost everywhere.

\begin{lemma}\label{l2}
For $f\ge0$ in $L^q$,
\begin{equation}\label{e23}
\HLS(f,f)\le\A(f)\le\HLS(f^\star,f^\star)
\le C_0\|f\|_q^2.
\end{equation}
In particular, $\A(f)$ is finite, and it equals $\HLS(f,f)$ for radial $f$.
\end{lemma}
\begin{proof}
This is the affine rearrangement inequality of Haddad and Ludwig
\cite{HL25b}. We recall its short kernel proof. For bounded,
compactly supported $f\ne0$, its autocorrelation is continuous,
compactly supported and positive near zero. Hence $S=S_\alpha f$
has a positive continuous radial function. Polar coordinates yield
\begin{equation}\label{e24}
\int_{\R^n}\int_{\R^n}\frac{f(x)f(y)}{\|x-y\|_K^{n-\alpha}}dxdy
=\int_{S^{n-1}}\rho_K^{n-\alpha}D_fd\xi.
\end{equation}
The kernel on the left rearranges to
$\big(\frac{|K|}{\omega_n}\big)^{\frac{n-\alpha}{n}}|x|^{-(n-\alpha)}$.
Taking $K=S$ and applying Riesz rearrangement yields
$n|S|\le\big(\frac{|S|}{\omega_n}\big)^{\frac{n-\alpha}{n}}\HLS(f^\star,f^\star)$,
which is the middle inequality. Apply this to
$f_j=\min(f,j)1_{B_j}\uparrow f$ as $j\rightarrow\infty$ and use monotone convergence
and $f_j^\star\le f^\star$ to treat general $f$.
Finally $\HLS(f,f)=\int_{S^{n-1}}D_fd\xi$. H\"older on the sphere proves the
first inequality, with equality when $D_f$ is constant. The last
inequality is classical HLS.
\end{proof}

\begin{proposition}
The two variational constants coincide:
$0<C_{\eps,\rm E}=C_{\eps,\rm A}=C_\eps\le C_0$.
For $f\ge0$ in $X_\eps$,
$\HLS(f,f)\le C_\eps\Scale{f}^2$ and
$\A(f)\le C_\eps\Scale{f}^2$.
\end{proposition}
\begin{proof}
Normalize any nonzero bounded compactly supported test to obtain
positivity. By normalization, equimeasurability and \eqref{e23},
\[
\HLS(f,f)\le\A(f)\le\HLS(f^\star,f^\star)
\le C_{\eps,\rm E}\Scale{f}^2\le C_0\Scale{f}^2.
\]
Taking suprema proves the assertions.
\end{proof}

\begin{lemma}\label{l3}
If $1<a<q<b$, every bounded sequence of nonnegative radial
nonincreasing functions in $L^a\cap L^b$ has a subsequence
converging almost everywhere and strongly in $L^q$.
\end{lemma}
\begin{proof}
Writing $L_p=\sup_j\int_{\R^n} f_j^p dx$ for $p=a,b$, monotonicity yields
\begin{equation}\label{e25}
f_j(r)\le\Big(\frac{L_a}{\omega_nr^n}\Big)^{\frac{1}{a}}.
\end{equation}
Helly selection on compact radial intervals yields an almost-everywhere
limit $f\in L^a\cap L^b$. Moreover,
\[
\int_{B_r}f_j^q dx\le L_b^{\frac{q}{b}}(\omega_nr^n)^{1-\frac{q}{b}},\quad
\int_{B_R^c}f_j^q dx\le L_a\Big(\frac{L_a}{\omega_nR^n}\Big)^{\frac{q-a}{a}}.
\]
The same bounds hold for $f$. They control the origin and infinity
uniformly, while \eqref{e25} yields dominated convergence on annuli.
Thus $f_j\rightarrow f$ in $L^q$ as $j\rightarrow\infty$.
\end{proof}

\begin{theorem}\label{t3}
There is $u_\eps=u_\eps^\star\in X_\eps$ with
$\M(u_\eps)=1$ and $\HLS(u_\eps,u_\eps)=\A(u_\eps)=C_\eps$.
Every normalized radial nonincreasing Euclidean maximizing sequence
has a subsequence converging to an extremal in $X_\eps$.
\end{theorem}
\begin{proof}
Rearrangement supplies such a maximizing sequence. Lemma~\ref{l3}
and HLS yield, after extraction, $f_j\rightarrow u$ almost everywhere and in
$L^q$ as $j\rightarrow\infty$, with $\HLS(u,u)=C_\eps>0$. Fatou yields $\M(u)\le1$.
Strict inequality would make $\frac{u}{\Scale{u}}$ a better competitor.
Thus $\M(u)=1$, and radiality yields affine attainment. For $p=a,b$,
the Br\'ezis-Lieb lemma \cite{BL} yields
$\int_{\R^n}|f_j-u|^p dx=\int_{\R^n} f_j^p dx-\int_{\R^n} u^p dx+o(1)$.
Adding the two identities and using normalization proves convergence
in $X_\eps$.
\end{proof}

\begin{remark}
The constraint and affine energy are invariant under translations
and $\SL(n)$: indeed $S_\alpha(f\circ T)=T^{-1}S_\alpha f$ for
$\det T=1$. Hence $u_\eps(T(\cdot-x_0))$ is also an affine extremal.
The converse is proved in Theorem~\ref{t8}.
\end{remark}

For the critical dilation $f_\lambda=\mathcal{D}_\lambda f$, both
energies are unchanged and, writing $A_f=\int_{\R^n}|f|^a dx$, $B_f=\int_{\R^n}|f|^b dx$,
\[
A_{f_\lambda}=\lambda^{-\frac{n\eps}{q}}A_f,\quad
B_{f_\lambda}=\lambda^{\frac{n\eps}{q}}B_f.
\]
For the affine energy, this follows directly from $D_{f_\lambda}=D_f$.

\begin{proposition}\label{p2}
For $0\ne f\in X_\eps$, set $s_f=(A_fB_f)^{\frac{1}{2q}}$.
Then $\min_{\lambda>0}\Scale{f_\lambda}=s_f$, attained uniquely at
\begin{equation}\label{e32}
\lambda_f=\Big(\frac{A_f}{B_f}\Big)^{\frac{q}{2n\eps}}\big(A_fB_f\big)^{\frac{1}{2n}}.
\end{equation}
At this scale $\int_{\R^n}\Big|\frac{f_{\lambda_f}}{s_f}\Big|^a dx
=\int_{\R^n}\Big|\frac{f_{\lambda_f}}{s_f}\Big|^b dx=1$.
\end{proposition}
\begin{proof}
For $s=\Scale{f_\lambda}$, arithmetic-geometric mean yields
\[
1=\frac{A_f\lambda^{-\frac{n\eps}{q}}}{2s^a}
 +\frac{B_f\lambda^{\frac{n\eps}{q}}}{2s^b}
\ge\frac{(A_fB_f)^{\frac{1}{2}}}{s^q}.
\]
Equality requires both summands before division by $2$ to equal $1$.
With $s=s_f$, this yields exactly \eqref{e32}.
\end{proof}

\begin{theorem}\label{t4}
Every nonnegative normalized extremal for either energy satisfies
\begin{equation*}
\int_{\R^n} u^a dx=\int_{\R^n} u^b dx=1.
\end{equation*}
\end{theorem}
\begin{proof}
Otherwise $A_u+B_u=2$ but $A_uB_u<1$. Proposition~\ref{p2}
would yield a dilation with normalization factor smaller than $1$
and hence a strictly larger normalized energy.
\end{proof}

\begin{proposition}
Every normalized nonnegative Euclidean extremal is positive almost
everywhere and satisfies
\begin{equation}\label{e35}
I_\alpha u:=\int_{\R^n}\frac{u(y)}{|x-y|^{n-\alpha}}dy
=\frac{C_\eps}{2q}(a u^{a-1}+b u^{b-1})\quad\text{a.e.}
\end{equation}
\end{proposition}
\begin{proof}
Both the constraint and energy are continuously Fr\'echet
differentiable on the real space $X_\eps$. Positivity of the kernel
makes $u$ a maximizer on the real constraint, and Theorem~\ref{t4} yields
$\M'(u)[u]=q$. Lagrange multipliers, tested against $u$, yield
$2\HLS(u,h)=(\frac{2C_\eps}{q})\M'(u)[h]$, hence \eqref{e35}.
Since $I_\alpha u>0$ wherever finite and
$I_\alpha u\in L^{q'}$ by HLS, while the right side vanishes at
$u=0$, positivity follows.
\end{proof}

\begin{theorem}\label{t5}
For $0<\eps<q-1$, $0<C_\eps<C_0$. The function
$\eps\mapsto C_\eps$ is strictly decreasing and
\begin{equation*}
C_\eps\rightarrow C_0\ \text{as }\eps\rightarrow0^+.
\end{equation*}
\end{theorem}
\begin{proof}
For $t>0$ with $t\ne1$, one has $\varphi_\eps(t)>t^q$.
A positive extremal $u\in L^q(\R^n)$ differs from $1$ on a set
of positive measure, so $\|u\|_q^q<1$ and $C_\eps<C_0$.
Since $\varphi_\eps(t)=t^q\cosh(\eps\log(t))$, an extremal for
$\eps_2>\eps_1$ belongs to $X_{\eps_1}$ and satisfies
$M_{\eps_1}(u)<1$. Renormalization yields
$C_{\eps_1}>C_{\eps_2}$. Finally, for the classical extremal $U$
in \eqref{e10}, dominated convergence yields
$M_\eps(\frac{U}{s})\rightarrow s^{-q}$ as $\eps\rightarrow0^+$ for each $s>0$, hence
$\Scale{U}\rightarrow1$ as $\eps\rightarrow0^+$. Thus $C_\eps\ge \frac{C_0}{\Scale{U}^2}\rightarrow C_0$ as $\eps\rightarrow0^+$.
\end{proof}

\begin{proposition}
The product formulations are
\begin{align}
C_\eps&=\sup_{0\ne f\ge0,\ f\in X_\eps}
\frac{\HLS(f,f)}{(A_fB_f)^{\frac{1}{q}}},\label{e38}\\
C_\eps&=\sup_{0\ne f\ge0,\ f\in X_\eps}
\frac{\A(f)}{(A_fB_f)^{\frac{1}{q}}}.\notag
\end{align}
For $f\ge0$, $\M(f)=1$, set
$\Delta_\eps(f)=C_\eps-\A(f)$ and
$\Delta_\eps^{\rm int}(f)=C_\eps(A_fB_f)^{\frac{1}{q}}-\A(f)\ge0$.
Then
\begin{gather}
\Delta_\eps(f)=\Delta_\eps^{\rm int}(f)
+C_\eps\big[1-(A_fB_f)^{\frac{1}{q}}\big],\notag\\
\Delta_\eps(f)\ge\frac{C_\eps}{4q}(A_f-B_f)^2.\label{e40}
\end{gather}
The same statements hold with $\HLS(f,f)$ in place of $\A(f)$.
\end{proposition}
\begin{proof}
Proposition~\ref{p2} minimizes the squared normalization factor at
$(A_fB_f)^{\frac{1}{q}}$ without changing either energy, proving the two
supremum formulas. The decomposition follows from the definitions of
$\Delta_\eps(f)$ and $\Delta_\eps^{\rm int}(f)$. Finally,
$A_f+B_f=2$ implies $A_fB_f=1-\frac{(A_f-B_f)^2}{4}$. Apply
$(1-t)^{\frac{1}{q}}\le1-\frac{t}{q}$.
\end{proof}

The Euclidean existence assertion is also covered by
Greco-Huang-Liu-Moroz~\cite[Theorem 1.2 and Remark 1.2]{GHLM}:
take $P=\frac{2}{a}$, $R=\frac{2b}{a}$, $\theta=\frac{a}{2q}$ and $f=v^P$.
Their parameter conditions become $a<q<b$, and their expression equals
\[
\frac{\Gamma(\frac{n-\alpha}{2})}{2^\alpha\pi^{\frac{n}{2}}\Gamma(\frac{\alpha}{2})}
\frac{\HLS(f,f)}{(A_fB_f)^{\frac{1}{q}}}.
\]
The proof above also records radial compactness and the identities
in Theorem~\ref{t4}, which are used below.

\section{Sharp asymptotics and scale selection}\label{s3}

Set
\[
\mathcal{R}_\eps=
\frac{\Gamma(n)^2\Gamma(\frac{n}{2}-\gamma\eps)
\Gamma(\frac{n}{2}+\gamma\eps)}{\Gamma(\frac{n}{2})^2\Gamma(n-\gamma\eps)\Gamma(n+\gamma\eps)}.
\]
All Gamma arguments are positive under \eqref{e1}.

\begin{proposition}
One has
\begin{gather}
C_0\mathcal{R}_\eps^{-\frac{1}{q}}\le C_\eps<C_0,\label{e42}\\
\limsup_{\eps\rightarrow0^+}\frac{C_0-C_\eps}{C_0\eps^2}
\le\kappa_{n,\alpha}.\label{e43}
\end{gather}
\end{proposition}
\begin{proof}
The beta integral
$\int_{\R^n}(1+|x|^2)^{-\sigma} dx
=\frac{\pi^{\frac{n}{2}}\Gamma(\sigma-\frac{n}{2})}{\Gamma(\sigma)}$ yields
$A_UB_U=\mathcal{R}_\eps$. Testing \eqref{e38} with $U$ proves
the lower bound in \eqref{e42}. The upper bound is Theorem~\ref{t5}.
Taylor expansion yields
\[
\log(\mathcal{R}_\eps)=\mathcal{V}_{n,\alpha}\eps^2+O(\eps^4),
\quad
\mathcal{R}_\eps^{-\frac{1}{q}}=1-\kappa_{n,\alpha}\eps^2+O(\eps^4)\ \text{as }\eps\rightarrow0^+,
\]
which proves \eqref{e43}. The coefficient is positive because
$\psi_1(t)=\sum_{k\ge0}(t+k)^{-2}$ is strictly decreasing.
\end{proof}

For the classical extremal family
$\mathcal{B}=\{c(\mathcal{D}_\rho U)(\cdot-y):c\ge0,\rho>0,y\in\R^n\}$,
the HLS stability theorem \cite[Theorem 1.1]{CLT24}, with
$s=\frac{\alpha}{2}$ and the kernel normalization absorbed into
$\beta_{n,\alpha}>0$, yields
\begin{equation}\label{e44}
C_0\|f\|_q^2-\HLS(f,f)
\ge\beta_{n,\alpha}\inf_{V\in\mathcal{B}}\|f-V\|_q^2,
\quad 0\le f\in L^q.
\end{equation}
Negative amplitudes cannot improve the distance for $f\ge0$, since
zero already belongs to $\mathcal{B}$.

\begin{lemma}\label{l4}
If $f=f^\star$ and $\|f\|_q=1$, then
\[
\inf_{\rho>0}\|f-\mathcal{D}_\rho U\|_q^2
\le\frac{4}{\beta_{n,\alpha}}(C_0-\HLS(f,f)).
\]
\end{lemma}
\begin{proof}
Rearrangement contraction in $L^q$ \cite{LL01} yields
$\|f-c\mathcal{D}_\rho U\|_q\le\|f-V\|_q$ for
$V=c(\mathcal{D}_\rho U)(\cdot-y)$.
Since $|c-1|\le\|f-V\|_q$, we have
$\|f-\mathcal{D}_\rho U\|_q\le2\|f-V\|_q$.
Take infima in \eqref{e44}, allowing $c\rightarrow0^+$.
\end{proof}

For $f\ge0$ in $X_\eps$ with $\|f\|_q=1$, define
\begin{equation*}
\mathcal{P}_\eps(f)=\bigg(\int_{\R^n} f^{q-\eps} dx\bigg)
                     \bigg(\int_{\R^n} f^{q+\eps} dx\bigg),
\end{equation*}
and let $\mu(f)$ and $\mathcal{V}(f)$ be the mean and variance of
$\log(f)$ under the probability measure $f^qdx$:
\[
\mu(f)=\int_{\R^n} f^q\log(f) dx,\quad
\mathcal{V}(f)=\int_{\R^n} f^q(\log(f))^2 dx-\mu(f)^2.
\]
The integrands are zero at $f=0$ and integrable, since
$t^q|\log(t)|^2\le C_\eps(t^{q-\eps}+t^{q+\eps})$.

\begin{lemma}
Writing $\Omega_f=\{f>0\}$, one has
\begin{equation}\label{e48}
\mathcal{P}_\eps(f)-1
=\iint_{\Omega_f^2}f(x)^qf(y)^q
 \big[\cosh(\eps(\log(f(x))-\log(f(y))))-1\big]dxdy,
\end{equation}
and consequently
\begin{equation}\label{e49}
\mathcal{P}_\eps(f)-1\ge\eps^2\mathcal{V}(f)\ge0.
\end{equation}
dilations satisfy
\begin{align*}
\mathcal{P}_\eps(\mathcal{D}_\rho f)&=\mathcal{P}_\eps(f),\notag\\
\mu(\mathcal{D}_\rho f)&=\mu(f)+\frac{n}{q}\log(\rho),\\
\mathcal{V}(\mathcal{D}_\rho f)&=\mathcal{V}(f).
\end{align*}
\end{lemma}
\begin{proof}
Symmetrizing the product integral proves \eqref{e48}.
The inequality $\cosh(t)-1\ge \frac{t^2}{2}$ and the identity
\[
\mathcal{V}(f)=\frac{1}{2}\iint_{\Omega_f^2}
 f(x)^qf(y)^q(\log(f(x))-\log(f(y)))^2dxdy
\]
yield \eqref{e49}. A change of variables proves the dilation identities:
the probability measure is transported and $\log(f)$ gains the
constant $(\frac{n}{q})\log(\rho)$.
\end{proof}

\begin{lemma}\label{l6}
Suppose $\eps_j\rightarrow0^+$ as $j\rightarrow\infty$, $0\le f_j\in X_{\eps_j}$,
$\|f_j\|_q=1$, and $f_j\rightarrow U$ strongly in $L^q$ and almost everywhere as $j\rightarrow\infty$.
Then
\begin{align}
\liminf_{j\rightarrow\infty}\frac{\mathcal{P}_{\eps_j}(f_j)-1}{\eps_j^2}
&\ge\mathcal{V}(U),\notag\\
\liminf_{j\rightarrow\infty}\mathcal{V}(f_j)&\ge\mathcal{V}(U).\label{e51-variance}
\end{align}
\end{lemma}
\begin{proof}
Since $U>0$, for almost every $(x,y)$ both values of $f_j$ are
eventually positive, and
\[
\frac{\cosh(\eps_j(\log(f_j(x))-\log(f_j(y))))-1}{\eps_j^2}
\rightarrow\frac{1}{2}(\log(U(x))-\log(U(y)))^2\ \text{as }j\rightarrow\infty.
\]
Extend the nonnegative integrands by zero outside $\Omega_{f_j}^2$.
Fatou's lemma applied to \eqref{e48} and to the variance integral
proves both assertions.
\end{proof}

\begin{lemma}\label{l7}
One has $\mathcal{V}(U)=\mathcal{V}_{n,\alpha}$.
\end{lemma}
\begin{proof}
The beta integral yields, for small $|s|$,
\begin{equation}\label{e52}
\int_{\R^n} U^{q+s} dx
=Z_n^{-1-\frac{s}{q}}\pi^{\frac{n}{2}}
  \frac{\Gamma(\frac{n}{2}+\gamma s)}{\Gamma(n+\gamma s)}.
\end{equation}
Twice differentiating its logarithm at zero yields the assertion.
Differentiation under the integral is justified by the explicit decay
of $U$ and integrability at exponents on both sides of $q$.
\end{proof}

The logarithmic integral estimates above determine the constant
and norm expansions in Theorem~\ref{t1}.
\begin{proof}
Let $\eps_j\rightarrow0^+$ as $j\rightarrow\infty$, $u_j\in\mathcal{E}_{\eps_j}$,
$s_j=\|u_j\|_q$ and $f_j=\frac{u_j}{s_j}$. Theorem~\ref{t4} yields
\begin{equation}\label{e53}
\mathcal{P}_{\eps_j}(f_j)=s_j^{-2q},\quad
C_{\eps_j}=\frac{\HLS(f_j,f_j)}{\mathcal{P}_{\eps_j}(f_j)^{\frac{1}{q}}}.
\end{equation}
H\"older, HLS and \eqref{e42} imply
\begin{equation}\label{e54}
1\le\mathcal{P}_{\eps_j}(f_j)
\le\Big(\frac{C_0}{C_{\eps_j}}\Big)^q\le\mathcal{R}_{\eps_j}
=1+\mathcal{V}_{n,\alpha}\eps_j^2+O(\eps_j^4)\ \text{as }j\rightarrow\infty.
\end{equation}
Also $C_{\eps_j}\le\HLS(f_j,f_j)\le C_0$, so the energy tends to
$C_0$. Lemma~\ref{l4} supplies $\rho_j>0$ such that
$g_j=\mathcal{D}_{\frac{1}{\rho_j}}f_j\rightarrow U$ in $L^q$ as $j\rightarrow\infty$.
For a subsequence realizing the liminf, extract almost everywhere
convergence and apply Lemmas~\ref{l6}-\ref{l7} and dilation invariance.
Together with \eqref{e54}, this yields
\begin{equation}\label{e55}
\mathcal{P}_{\eps_j}(f_j)
=1+\mathcal{V}_{n,\alpha}\eps_j^2+o(\eps_j^2)\ \text{as }j\rightarrow\infty.
\end{equation}
Since $\frac{C_{\eps_j}}{C_0}\le\mathcal{P}_{\eps_j}(f_j)^{-\frac{1}{q}}$,
\eqref{e55} and \eqref{e43} prove \eqref{e13}.
Equation \eqref{e53} yields
$s_j^2=1-\kappa_{n,\alpha}\eps_j^2+o(\eps_j^2)$, proving
\eqref{e14}. The initial sequence was arbitrary.
\end{proof}

\begin{lemma}\label{l8}
Suppose $\eps_j\rightarrow0^+$ as $j\rightarrow\infty$, $0\le g_j\in X_{\eps_j}$,
$\|g_j\|_q=1$, $g_j\rightarrow U$ in $L^q$ and almost everywhere as $j\rightarrow\infty$, and
$\frac{\mathcal{P}_{\eps_j}(g_j)-1}{\eps_j^2}\rightarrow\mathcal{V}_{n,\alpha}$ as $j\rightarrow\infty$.
Then $\mu(g_j)\rightarrow\mu(U)$ as $j\rightarrow\infty$, $\mathcal{V}(g_j)\rightarrow\mathcal{V}_{n,\alpha}$ as $j\rightarrow\infty$,
and, for $L_j^\pm=\int_{\R^n} g_j^{q\pm\eps_j} dx$,
\begin{equation}\label{e57}
\log(L_j^\pm)=\pm\eps_j\mu(g_j)+O(\eps_j^2)\ \text{as }j\rightarrow\infty.
\end{equation}
\end{lemma}
\begin{proof}
Equations \eqref{e49} and \eqref{e51-variance} yield convergence
of the variances. Choose a bounded annulus $E$ on which
$0<c\le U\le C$. On
$E_j=E\cap\{\frac{c}{2}\le g_j\le2C\}$, dominated convergence yields
$\int_{E_j}g_j^q dx\rightarrow\int_E U^q dx>0$ as $j\rightarrow\infty$.
Thus $\int_{\{|\log(g_j)|\le M\}}g_j^qdx\ge\eta>0$ for some fixed
$M,\eta$ and all large $j$, whence
$\mathcal{V}(g_j)\ge\eta(|\mu(g_j)|-M)_+^2$.
The means are bounded, and therefore
\begin{equation}\label{e58}
\int_{\R^n} g_j^q(\log(g_j))^2 dx
=\mathcal{V}(g_j)+\mu(g_j)^2\le K.
\end{equation}
For the truncation $T_L(t)=\max(-L,\min(t,L))$, strong $L^1$
convergence of $g_j^q$ and almost everywhere convergence yield
$\int_{\R^n} g_j^qT_L(\log(g_j)) dx\rightarrow\int_{\R^n} U^qT_L(\log(U)) dx$ as $j\rightarrow\infty$.
Equation \eqref{e58} bounds each truncation remainder by $\frac{K}{L}$
(with an analogous bound for $U$), proving convergence of the means.
Finally Jensen yields $\log(L_j^\pm)\mp\eps_j\mu(g_j)\ge0$.
Their sum is $\log(\mathcal{P}_{\eps_j}(g_j))=O(\eps_j^2)$,
which proves \eqref{e57}.
\end{proof}

\begin{lemma}
Every choice $u_\eps\in\mathcal{E}_\eps$ satisfies
\begin{align}
\int_{\R^n} u_\eps^q\log(u_\eps)dx&\rightarrow0,\label{e59}\\
\int_{\R^n} u_\eps^q(\log(u_\eps))^2 dx&\rightarrow\mathcal{V}_{n,\alpha}
\ \text{as }\eps\rightarrow0^+.\label{e59-second-moment}
\end{align}
The scales $\rho_j$ chosen in the proof of \eqref{e13} and
\eqref{e14} satisfy $\rho_j\rightarrow\rho_*$ as $j\rightarrow\infty$.
\end{lemma}
\begin{proof}
Use the notation of that proof and extract almost everywhere
convergence $g_j\rightarrow U$ as $j\rightarrow\infty$. Lemma~\ref{l8} applies by \eqref{e55}.
Theorem~\ref{t4} yields
$L_j^-=\rho_j^{n\eps_j/q}s_j^{-(q-\eps_j)}$ and
$L_j^+=\rho_j^{-n\eps_j/q}s_j^{-(q+\eps_j)}$.
Taking logarithmic differences in \eqref{e57},
\begin{equation}\label{e60}
\frac{n}{q}\log(\rho_j)+\log(s_j)=-\mu(g_j)+O(\eps_j)\ \text{as }j\rightarrow\infty.
\end{equation}
Since $s_j\rightarrow1$ as $j\rightarrow\infty$ and $\mu(g_j)\rightarrow\mu(U)$ as $j\rightarrow\infty$,
$\rho_j\rightarrow\exp\big[-\frac{q\mu(U)}{n}\big]$ as $j\rightarrow\infty$.
Differentiating \eqref{e52} once yields
\[
\mu(U)=-\frac{1}{q}\log(Z_n)
       -\gamma\Big(\psi(n)-\psi\Big(\frac{n}{2}\Big)\Big),
\]
so $q\gamma=n$ identifies this limit with \eqref{e11}.
The strong continuity and isometry of critical dilation yield
\[
\|u_j-U_*\|_q\le s_j\|g_j-U\|_q+|s_j-1|
 +\|\mathcal{D}_{\rho_j}U-\mathcal{D}_{\rho_*}U\|_q\rightarrow0\ \text{as }j\rightarrow\infty.
\]
Under $u_j=s_j\mathcal{D}_{\rho_j}g_j$, the logarithm acquires the
shift $\log(s_j)+(\frac{n}{q})\log(\rho_j)$. By \eqref{e60} its mean becomes
$O(\eps_j)$ and its variance remains $\mathcal{V}(g_j)$.
Multiplication by $s_j^q\rightarrow1$ as $j\rightarrow\infty$ proves \eqref{e59} and
\eqref{e59-second-moment}. Subsequence extraction from any initial
sequence proves the limits and completes \eqref{e15}.
\end{proof}

Since Theorem~\ref{t1} applies to every choice of radial extremals,
\begin{equation}\label{e61}
\sup_{u\in\mathcal{E}_\eps}\|u-U_*\|_q\rightarrow0
\ \text{as }\eps\rightarrow0^+.
\end{equation}
Otherwise a violating sequence would contradict the theorem.

\begin{corollary}
For any normalized affine extremals $w_\eps\ge0$,
$w_\eps^\star\rightarrow U_*$ strongly in $L^q$ as $\eps\rightarrow0^+$.
\end{corollary}
\begin{proof}
The rearrangement chain yields
$C_\eps=\A(w_\eps)\le\HLS(w_\eps^\star,w_\eps^\star)\le C_\eps$.
Thus $w_\eps^\star\in\mathcal{E}_\eps$, and Theorem~\ref{t1} applies.
\end{proof}

\section{Extremals and compactness}\label{s4}

We first obtain Euclidean compactness, then use a quantitative
rearrangement estimate to control the affine freedom. Write $B=B_1$.
For a measurable star set $K$ with positive finite radial function
and $|K|=\omega_n$, define
\[
k_K(z)=\|z\|_K^{-(n-\alpha)},\quad
\mathcal{H}_K(f,h)=\int_{\R^n}\int_{\R^n} f(x)h(y)k_K(x-y)dxdy.
\]
Since $k_K^\star(z)=|z|^{-(n-\alpha)}$, Riesz rearrangement yields
\begin{equation}\label{e62}
|\mathcal{H}_K(f,h)|\le C_0\|f\|_q\|h\|_q.
\end{equation}
For $0\ne f\ge0$ in $L^q$, set
\begin{equation}\label{e63}
K_f=\bigg(\frac{\omega_n}{|S_\alpha f|}\bigg)^{\frac{1}{n}}S_\alpha f.
\end{equation}
This is an origin-symmetric star set of volume $\omega_n$.
Indeed, $f\ge c1_E$ for some bounded positive-measure set $E$.
Continuity at zero of the autocorrelation of $1_E$ yields a uniform
positive lower bound for $D_f$. Finiteness follows from
Lemma~\ref{l2}, and symmetry from the evenness of autocorrelation.

\begin{proposition}
For $f\ge0$ in $L^q$,
\begin{equation*}
\A(f)=\sup_{|K|=\omega_n}\mathcal{H}_K(f,f),
\end{equation*}
with attainment at $K_f$ when $f\ne0$. For $f,h\ge0$,
\begin{equation}\label{e65}
|\A(f)-\A(h)|\le C_0(\|f\|_q+\|h\|_q)\|f-h\|_q.
\end{equation}
\end{proposition}
\begin{proof}
Identity \eqref{e24} remains valid for measurable star sets by
Tonelli. H\"older and $\int_{S^{n-1}}\rho_K^n d\xi=n\omega_n$ yield
$\mathcal{H}_K(f,f)\le\A(f)$, with equality for \eqref{e63}.
For each $K$, \eqref{e62} bounds
$|\mathcal{H}_K(f,f)-\mathcal{H}_K(h,h)|$ by the right side of
\eqref{e65}. Take suprema.
\end{proof}

\begin{lemma}\label{l10}
For $|K|=\omega_n$ and $f\ge0$ in $L^q$,
\[
|\mathcal{H}_K(f,f)-\HLS(f,f)|
\le C_0\|f\|_q^2
\bigg(\frac{|K\triangle B|}{\omega_n}\bigg)^{\frac{n-\alpha}{n}}.
\]
\end{lemma}
\begin{proof}
For $G=|k_K-k_B|$, homogeneity yields
\[
|\{G>t\}|=t^{-\frac{n}{n-\alpha}}\Theta_K,\quad
\Theta_K=\frac{1}{n}\int_{S^{n-1}}
|\rho_K^{n-\alpha}-1|^{\frac{n}{n-\alpha}}d\xi
\le|K\triangle B|.
\]
The last inequality follows from $|x^\beta-1|\le|x-1|^\beta$
for $0<\beta<1$. Thus
$G^\star(z)=(\frac{\Theta_K}{\omega_n})^{\frac{n-\alpha}{n}}|z|^{-(n-\alpha)}$.
Riesz rearrangement followed by HLS proves the claim.
\end{proof}

\begin{theorem}\label{t6}
For fixed $0<\eps<q-1$, every normalized nonnegative Euclidean
maximizing sequence has a subsequence converging strongly in
$X_\eps$ to an extremal after translations.
\end{theorem}
\begin{proof}
By the Euclidean version of \eqref{e40} and normalization,
\begin{equation}\label{e67}
 \int_{\R^n} f_j^a dx\rightarrow1,\quad
 \int_{\R^n} f_j^b dx\rightarrow1
 \quad\text{as }j\rightarrow\infty.
\end{equation}
Apply Lions' concentration compactness alternatives \cite{L84,L85}
to $d\nu_j=\frac{1}{2}(f_j^a+f_j^b) dx$.
Choose $p_+\in(q,\min\{b,2\})$ and set
$r_p=\frac{p}{2(p-1)}$ for $1<p<2$.
Young's inequality bounds the bilinear form with kernel $k$ by
$\|g\|_p\|h\|_p\|k\|_{r_p}$. Since
$r_{p_+}<\frac{n}{n-\alpha}<r_a$, the uniform $L^{p_+}$ and $L^a$
bounds imply
\begin{align}
\sup_j\int_{\R^n}\int_{|x-y|<s}
\frac{f_j(x)f_j(y)}{|x-y|^{n-\alpha}}dxdy&\rightarrow0
\ \text{as }s\rightarrow0^+,\label{e68}\\
\sup_j\int_{\R^n}\int_{|x-y|>R}
\frac{f_j(x)f_j(y)}{|x-y|^{n-\alpha}}dxdy&\rightarrow0
\ \text{as }R\rightarrow\infty.\label{e69}
\end{align}
These bounds also apply to functions dominated by $f_j$.

In the vanishing case, partition space into fixed cubes $Q_i$.
For fixed $0<s<R$, only finitely many neighboring cubes interact
under the truncated kernel. With $c_{j,i}=\int_{Q_i}f_j^a dx$,
H\"older and the inequality $2xy\le x^2+y^2$ bound this energy by
\[
C_{s,R}\sum_i c_{j,i}^{\frac{2}{a}}
\le C_{s,R}(\sup_i c_{j,i})^{\frac{2}{a}-1}\sum_i c_{j,i}\rightarrow0\ \text{as }j\rightarrow\infty.
\]
Here vanishing yields $\sup_i c_{j,i}\rightarrow0$ as $j\rightarrow\infty$ and $a<2$.
Together with \eqref{e68}-\eqref{e69}, this would force the entire
energy to zero.

In the dichotomy case, write $f_j=g_j+h_j+z_j$ as restrictions
to disjoint sets, with the supports of $g_j,h_j$ separating,
$\|z_j\|_{X_\eps}\rightarrow0$ as $j\rightarrow\infty$, and limiting $\nu_j$ masses
$\vartheta,1-\vartheta\in(0,1)$. After extraction write
$\int_{\R^n} g_j^a dx\rightarrow A_1$ as $j\rightarrow\infty$, $\int_{\R^n} h_j^a dx\rightarrow A_2$ as $j\rightarrow\infty$ and similarly
$\int_{\R^n} g_j^b dx\rightarrow B_1$ as $j\rightarrow\infty$, $\int_{\R^n} h_j^b dx\rightarrow B_2$ as $j\rightarrow\infty$. Then
\[
A_1+A_2=B_1+B_2=1,\quad
A_1+B_1=2\vartheta,\quad A_2+B_2=2(1-\vartheta).
\]
HLS removes all terms involving $z_j$, and \eqref{e69} removes
the $g_j,h_j$ interaction. The product bound \eqref{e38} then yields
\[
1\le\sum_{i=1}^2(A_iB_i)^{\frac{1}{q}}
\le\bigg(\sum_iA_i^{\frac{2}{q}}\bigg)^{\frac{1}{2}}
   \bigg(\sum_iB_i^{\frac{2}{q}}\bigg)^{\frac{1}{2}}\le1.
\]
Equality throughout would require, since $\frac{2}{q}>1$, both pairs to be
$(1,0)$ or $(0,1)$, and equality in Cauchy-Schwarz would require
the same nonzero index. This contradicts $0<\vartheta<1$.

Thus, after translations, $\nu_j$ is tight. Extract a common weak
limit $v$ in $L^a\cap L^b$, hence in $L^q$. The bound
\begin{equation}\label{e71}
\int_F f_j^q dx\le
\bigg(\int_F f_j^a dx\int_F f_j^b dx\bigg)^{\frac{1}{2}}\le\nu_j(F)
\end{equation}
yields $L^q$ tightness, also for $v$ by weak lower semicontinuity.
We claim that
$\HLS(f_j,f_j)\rightarrow\HLS(v,v)$ as $j\rightarrow\infty$. Weak lower semicontinuity bounds the
$L^a$ and $L^{p_+}$ norms of $v$, so the same Young estimates yield
the near and far bounds for $v$. Remove those parts of
the kernel, then approximate the remaining annular kernel by a
continuous compactly supported one using Young's inequality.
The kernels can be bounded by a fixed multiple of the Riesz kernel,
so \eqref{e71} and HLS allow restriction of both variables to a
fixed ball with uniformly small error. On its product, uniform
approximation by finite sums $\phi(x)\psi(y)$ reduces convergence
to weak $L^q$ convergence. Uniform local $L^1$ bounds control the
approximation remainder, proving the claim.

Consequently $\HLS(v,v)=C_\eps$. Weak lower semicontinuity yields
$\M(v)\le1$, and normalization of $v$ rules out strict inequality.
Theorem~\ref{t4} yields $\int_{\R^n} v^a dx=\int_{\R^n} v^b dx=1$. Equation \eqref{e67} and uniform
convexity of $L^a,L^b$ now yield strong convergence in $X_\eps$.
\end{proof}

For a set $K$ of volume $\omega_n$, define its ellipsoid asymmetry by
\begin{equation*}
\mathfrak{a}(K)=\frac{1}{\omega_n}
\inf_{\substack{E\text{ centered ellipsoid}\\|E|=\omega_n}}
|K\triangle E|
=\frac{1}{\omega_n}\inf_{T\in\SL(n)}|TK\triangle B|\le2.
\end{equation*}
For positive finite-measure sets, set
$\mathcal{T}(E,F,G)=\int_{\R^n}\int_{\R^n}1_E(x)1_F(y)1_G(x-y)dxdy$ and
\[
\mathfrak{D}(E,F,G)=\inf_{\substack{T\in\SL(n)\\y,z\in\R^n}}
\max\{|E\triangle(TE^\star+y)|,
|F\triangle(TF^\star+z)|,|G\triangle(TG^\star+y-z)|\}.
\]
Let $r_i$ be their volume radii. If
$r_k\le(1-\tau)(r_i+r_j)$ for every permutation and
$\min_i r_i\ge\tau\max_i r_i$, Christ's theorem
\cite[Theorem 1]{Chr}, applied to $(E,-F,G)$, yields
\begin{equation}\label{e73}
\mathcal{T}(E^\star,F^\star,G^\star)-\mathcal{T}(E,F,G)
\ge c_{n,\tau}\mathfrak{D}(E,F,G)^2.
\end{equation}

\begin{lemma}\label{l11}
If $K=-K$ and $|K|=\omega_n$, then
\begin{equation}\label{e74}
\inf_{\substack{E\text{ centered ellipsoid},\ |E|=\omega_n\\z\in\R^n}}
|K\triangle(z+E)|\ge\frac{\omega_n}{3}\mathfrak{a}(K).
\end{equation}
Moreover, $\mathfrak{a}(K)=0$ exactly when $K$ agrees almost everywhere
with a centered ellipsoid.
\end{lemma}
\begin{proof}
If $\eta=|K\triangle(z+E)|$, origin symmetry yields
$|(z+E)\triangle(-z+E)|\le2\eta$.
Mapping $E$ to a ball and slicing shows that its symmetric difference
from a translate increases with the displacement. Thus
$|E\triangle(z+E)|\le2\eta$ and $|K\triangle E|\le3\eta$.
If $\mathfrak{a}(K)=0$, choose $E_j\rightarrow K$ in measure as $j\rightarrow\infty$.
Their smallest semiaxes cannot tend to zero, since slicing would
yield $|E_j\cap B_R|\rightarrow0$ as $j\rightarrow\infty$ for each fixed $R$, contradicting $|K|>0$.
Their fixed volume also bounds the largest semiaxes. A subsequence
therefore converges to an ellipsoid equal to $K$ almost everywhere.
The converse is immediate.
\end{proof}

\begin{lemma}\label{l12}
For $0\ne u=u^\star\in L^q$, there are an $L^q$ neighborhood
$\mathcal{U}_u$ and $c_u>0$ such that $f\ge0$,
$f^\star\in\mathcal{U}_u$ and $K=-K$, $|K|=\omega_n$, imply
\begin{equation*}
\HLS(f^\star,f^\star)-\mathcal{H}_K(f,f)
\ge c_u\mathfrak{a}(K)^2.
\end{equation*}
\end{lemma}
\begin{proof}
Let $m_v(s)=|\{v>s\}|$. Choose a continuity point $s_0>0$ with
$0<m_u(s_0)<\infty$ and set $R_0=(\frac{m_u(s_0)}{\omega_n})^{\frac{1}{n}}$.
There are intervals $I\ni s_0$, $J\ni R_0$ of positive length
and a neighborhood $\mathcal{U}_u$ such that
\[
\bigg(\frac{m_v(s)}{\omega_n}\bigg)^{\frac{1}{n}},\ \bigg(\frac{m_v(t)}{\omega_n}\bigg)^{\frac{1}{n}},\ r
\in\bigg[\frac{4R_0}{5},\frac{6R_0}{5}\bigg]
\quad(s,t\in I,\ r\in J,\ v\in\mathcal{U}_u).
\]
Indeed choose the endpoints of $I$ at nearby continuity points of
$m_u$. Convergence of their level-set volumes under $L^q$
convergence, followed by monotonicity, yields uniform control on $I$.
These radii satisfy the hypotheses of \eqref{e73} with $\tau=\frac{1}{5}$.

Writing $E_s=\{f>s\}$ and using
$k_K(z)=(n-\alpha)\int_0^\infty r^{-(n-\alpha)-1}1_{rK}(z)dr$,
layer cake yields
\begin{equation}\label{e76}
\HLS(f^\star,f^\star)-\mathcal{H}_K(f,f)=(n-\alpha)\int_{(0,\infty)^3}r^{-(n-\alpha)-1}
[\mathcal{T}(E_s^\star,E_t^\star,rB)-\mathcal{T}(E_s,E_t,rK)]drdtds.
\end{equation}
The integrand is nonnegative by rearrangement and both energies are
finite by \eqref{e62}. On $I\times I\times J$, \eqref{e74} yields
$\mathfrak{D}(E_s,E_t,rK)\ge \frac{r^n\omega_n\mathfrak{a}(K)}{3}$.
Restricting \eqref{e76} to this box and applying \eqref{e73} proves
the result, with
$c_u=(n-\alpha)c_{n,\frac{1}{5}}\omega_n^2|I|^2
\frac{1}{9}\int_J r^{n+\alpha-1}dr>0$.
\end{proof}

\begin{theorem}\label{t7}
For fixed $0<\eps<q-1$, there is $\gamma_\eps>0$ such that
every nonnegative normalized $f$ satisfies
\begin{equation*}
\Delta_\eps(f)=C_\eps-\A(f)
\ge\gamma_\eps\mathfrak{a}(K_f)^2.
\end{equation*}
\end{theorem}
\begin{proof}
Compactness of $\mathcal{E}_\eps$ yields a finite cover by the
neighborhoods in Lemma~\ref{l12}, with a common constant $c_\eps>0$.
By Theorem~\ref{t3}, every normalized radial function with Euclidean
deficit below some $\delta_\eps>0$ lies in this cover.
The decomposition
\[
\Delta_\eps(f)=[C_\eps-\HLS(f^\star,f^\star)]
 +[\HLS(f^\star,f^\star)-\mathcal{H}_{K_f}(f,f)]
\]
has nonnegative terms. If $\Delta_\eps(f)<\delta_\eps$, apply
Lemma~\ref{l12}. Otherwise use $\mathfrak{a}(K_f)\le2$.
Thus $\gamma_\eps=\min(c_\eps,\frac{\delta_\eps}{4})$ works.
\end{proof}

\begin{theorem}\label{t8}
The normalized nonnegative affine extremals are exactly
$\mathcal{M}_{\eps,\rm A}$. Every normalized Euclidean extremal
is a translate of a member of $\mathcal{E}_\eps$.
\end{theorem}
\begin{proof}
A Euclidean extremal has equality in Riesz rearrangement. Its energy
is finite and positive, and the kernel is strictly radially decreasing.
The equality statement in \cite[Theorem 3.9]{LL01} therefore makes it
a translate of its radial rearrangement. For an affine extremal,
Theorem~\ref{t7} and Lemma~\ref{l11} show that $K_f$ is an ellipsoid.
Choose $L\in\SL(n)$ with $LK_f=B$ and set $g=f\circ L^{-1}$.
Then $\M(g)=1$ and
$\HLS(g,g)=\mathcal{H}_{K_f}(f,f)=\A(f)=C_\eps$,
so the Euclidean classification applies. Conversely, invariance
under translations and $\SL(n)$ proves attainment on every such orbit.
\end{proof}

\begin{theorem}\label{t9}
For fixed $0<\eps<q-1$, every normalized nonnegative affine
maximizing sequence admits a subsequence, $T_j\in\SL(n)$,
$x_j\in\R^n$, and $u\in\mathcal{E}_\eps$ such that
\begin{equation*}
f_j(T_j\cdot+x_j)\rightarrow u\quad\text{strongly in }X_\eps\ \text{as }j\rightarrow\infty.
\end{equation*}
\end{theorem}
\begin{proof}
Theorem~\ref{t7} yields $\mathfrak{a}(K_{f_j})\rightarrow0$ as $j\rightarrow\infty$.
Choose $L_j\in\SL(n)$ with $|L_jK_{f_j}\triangle B|\rightarrow0$ as $j\rightarrow\infty$ and set
$g_j=f_j\circ L_j^{-1}$. Then
$\mathcal{H}_{L_jK_{f_j}}(g_j,g_j)=\A(f_j)\rightarrow C_\eps$ as $j\rightarrow\infty$.
Lemma~\ref{l10} makes $g_j$ a Euclidean maximizing sequence.
Theorem~\ref{t6} yields strong convergence modulo translations,
and Theorem~\ref{t8} allows the limit to be centered and radial.
\end{proof}

\section{Euclidean stability}\label{s5}
The Euclidean version of \eqref{e40} yields two linearized
constraints along a sequence violating stability.
Their half sum removes the
constant mode, while their rescaled difference controls dilation.
We combine this observation with the classical HLS spectrum and a
convex-remainder normalization, following the local-to-global strategy
of \cite{BE,CFW}. See also \cite{Car,DEFFL,K23}.
Throughout this section $0<\eps<\min\{q-1,2-q\}$ is sufficiently small.

On the real Hilbert space $L^2(\R^n)$, write
$\langle z,w\rangle=\int_{\R^n}z(x)w(x) dx$.
We use $\rightharpoonup$ for weak convergence and $z\perp w$
for $\langle z,w\rangle=0$.
Set $s_\alpha=\frac{2n}{\alpha}=\frac{2q}{2-q}$, so $\frac{1}{q}=\frac{1}{2}+\frac{1}{s_\alpha}$.
For real $\psi\in L^{s_\alpha}$, define the selfadjoint operator
$\mathcal{K}_\psi$ on real $L^2(\R^n)$ by
\begin{equation*}
 \langle z,\mathcal{K}_\psi w\rangle=\HLS(\psi z,\psi w).
\end{equation*}

\begin{lemma}\label{l13}
The operator $\mathcal{K}_\psi$ is compact, and
\begin{equation}\label{e81}
 \|\mathcal{K}_\psi-\mathcal{K}_\eta\|_{L^2\rightarrow L^2}
 \le C_0(\|\psi\|_{s_\alpha}+\|\eta\|_{s_\alpha})
                  \|\psi-\eta\|_{s_\alpha}.
\end{equation}
\end{lemma}
\begin{proof}
H\"older and HLS yield boundedness and \eqref{e81}.
For bounded compactly supported $\psi$, removing $|x-y|<\delta$
changes the operator norm by at most
$\|\psi\|_\infty^2\int_{|z|<\delta}|z|^{-(n-\alpha)}dz=o(1)$,
by Young's inequality. The remaining kernel is square integrable.
Approximation in $L^{s_\alpha}$ and \eqref{e81} prove the assertion.
\end{proof}

For $u\in\mathcal{E}_\eps$, set
\[
 W_{\eps,u}=\frac{C_\eps}{2q}
 [a(a-1)u^{a-2}+b(b-1)u^{b-2}],\quad
 \psi_{\eps,u}=W_{\eps,u}^{-\frac{1}{2}}.
\]
Positivity of $u$ and the arithmetic-geometric mean inequality yield
$W_{\eps,u}\ge c_\eps u^{q-2}$, hence
$\psi_{\eps,u}\in L^{s_\alpha}$.
Write $W_*=C_0(q-1)U_*^{q-2}$ and $\psi_*=W_*^{-\frac{1}{2}}$.
The scale and logarithmic integral formulas of Section~\ref{s3} yield
\begin{align}
 \int_{\R^n} U_*^q dx=1,\quad \int_{\R^n} U_*^q\log(U_*) dx&=0,\label{e84}\\
 \int_{\R^n} U_*^q(\log(U_*))^2 dx&=\mathcal{V}_{n,\alpha}.
 \label{e84-second-moment}
\end{align}
For $1\le i\le n$, let $\partial_i$ denote differentiation
with respect to the $i$th coordinate. The centering tests
$\chi_i=U_*^{q-2}\partial_iU_*$ satisfy
$\chi_i,\partial_j\chi_i,\partial_j\partial_k\chi_i\in L^{q'}$.

\begin{lemma}\label{l14}
Suppose $\eps_j\rightarrow0^+$ as $j\rightarrow\infty$ and $u_j\in\mathcal{E}_{\eps_j}$.
Set $\psi_j=\psi_{\eps_j,u_j}$ and $a_j=q-\eps_j$, $b_j=q+\eps_j$.
Then
\begin{equation}\label{e86}
 \psi_j\rightarrow\psi_*\quad\text{in }L^{s_\alpha}\ \text{as }j\rightarrow\infty,
\end{equation}
and the following convergences hold in $L^2$:
\begin{gather}
 v_j^0=\frac{\psi_j}{2}(a_ju_j^{a_j-1}+b_ju_j^{b_j-1})
 \rightarrow v_*^0=\frac{q}{\sqrt{C_0(q-1)}}U_*^{\frac{q}{2}},\notag\\
 v_j^1=\frac{\psi_j}{2\eps_j}(b_ju_j^{b_j-1}-a_ju_j^{a_j-1})
 \rightarrow v_*^1=\frac{U_*^{\frac{q}{2}}(1+q\log(U_*))}{\sqrt{C_0(q-1)}},
 \label{e88}\\
 \psi_j\chi_i\rightarrow\psi_*\chi_i,\quad 1\le i\le n\ \text{as }j\rightarrow\infty.
 \label{e89}
\end{gather}
\end{lemma}
\begin{proof}
Theorem~\ref{t1} yields $u_j\rightarrow U_*$ in $L^q$ as $j\rightarrow\infty$. On an almost-everywhere
convergent subsequence, positivity of $U_*$ implies pointwise convergence
of these weights and vectors. Uniformly for small $\eps_j$,
$|\psi_j|^{s_\alpha}\le Cu_j^q$. Since $u_j^q\rightarrow U_*^q$ in $L^1$ as $j\rightarrow\infty$,
Vitali's theorem proves \eqref{e86}, including control at infinity.

For the divided difference, let
\[
 R_j^0=u_j^q\cosh(\eps_j\log(u_j)),\quad
 R_j^1=\frac{u_j^q}{\eps_j^2}\big(\cosh(\eps_j\log(u_j))-1\big).
\]
Theorem~\ref{t4} and \eqref{e14} yield $\int_{\R^n} R_j^0 dx=1$ and
$\int_{\R^n} R_j^1 dx\rightarrow\frac{\mathcal{V}_{n,\alpha}}{2}$ as $j\rightarrow\infty$. Their pointwise limits are
$U_*^q$ and $\frac{U_*^q(\log(U_*))^2}{2}$, respectively, with these same
integrals by \eqref{e84}-\eqref{e84-second-moment}.
Thus both nonnegative densities converge in $L^1$.
Writing $t=\eps_j\log(u_j)$, one has
\[
 v_j^0=\psi_ju_j^{q-1}(q\cosh(t)+\eps_j\sinh(t)),\quad
 v_j^1=\psi_ju_j^{q-1}\Big(\cosh(t)+\Big(\frac{q}{\eps_j}\Big)\sinh(t)\Big).
\]
The bounds $W_{\eps_j,u_j}\ge cu_j^{q-2}\cosh(t)$ and
$\frac{\sinh^2(t)}{\cosh(t)}\le2(\cosh(t)-1)$ yield
$|v_j^0|^2\le CR_j^0$ and $|v_j^1|^2\le C(R_j^0+R_j^1)$.
Vitali proves \eqref{e88}. H\"older and \eqref{e86} yield
\eqref{e89}. Applying this argument to a further subsequence of every
subsequence proves the stated sequence convergences.
\end{proof}

\begin{lemma}
There is $\sigma_*>0$ such that
\begin{equation}\label{e90}
 \|z\|_2^2-\HLS(\psi_*z,\psi_*z)\ge\sigma_*\|z\|_2^2
\end{equation}
whenever $z\perp v_*^0,v_*^1,\psi_*\chi_1,\ldots,\psi_*\chi_n$.
\end{lemma}
\begin{proof}
Let $S^n=\{\omega\in\R^{n+1}:|\omega|=1\}$.
The stereographic map
\[
 \mathcal{S}_{\rho_*}(x)=
 \bigg(\frac{2\rho_*x}{1+\rho_*^2|x|^2},
       \frac{1-\rho_*^2|x|^2}{1+\rho_*^2|x|^2}\bigg)
\]
sends $U_*^qdx$ to normalized surface measure on $S^n$.
Under the isometry $F\mapsto U_*^{\frac{q}{2}}(F\circ\mathcal{S}_{\rho_*})$,
the eigenvalues of $\mathcal{K}_{\psi_*}$ on degree-$\ell$ spherical
harmonics are, by the Funk-Hecke formula \cite{FL},
\[
 \lambda_\ell=\frac{1}{q-1}\prod_{j=0}^{\ell-1}
 \frac{n-\alpha+2j}{n+\alpha+2j}.
\]
The normalization follows also from $I_\alpha U_*=C_0U_*^{q-1}$.
Thus $\lambda_0=(q-1)^{-1}>1$, $\lambda_1=1$, and
$\lambda_\ell\le\lambda_2=\frac{n-\alpha+2}{n+\alpha+2}<1$ for $\ell\ge2$.
The constant mode is removed by $v_*^0$. The degree-one modes are
\[
 Z_i=\sqrt{W_*}\,\partial_iU_*,\quad
 Z_D=\sqrt{W_*}\,\Lambda U_*,\quad
 \Lambda F=\frac{n}{q}F+x\cdot\nabla F.
\]
Here $\nabla F=(\partial_1F,\ldots,\partial_nF)$ is the Euclidean
gradient. Since $\psi_*\chi_i$ is a nonzero multiple of $Z_i$, the translation
modes are removed. Dilation invariance and
$\mu(\mathcal{D}_{e^t}U_*)=\frac{nt}{q}$ yield
\[
 \int_{\R^n} U_*^{q-1}\Lambda U_* dx=0,\quad
 \int_{\R^n} U_*^{q-1}\log(U_*)\Lambda U_* dx=\frac{n}{q^2},
 \quad \langle v_*^1,Z_D\rangle=\frac{n}{q}\ne0.
\]
Write $z=ce_D+z_\perp$, where $e_D=\frac{Z_D}{\|Z_D\|_2}$ and
$z_\perp$ has only degrees at least two. Orthogonality to $v_*^1$
implies $|c|\le\frac{\|v_*^1\|_2\|z_\perp\|_2}{|\langle v_*^1,e_D\rangle|}$.
The gap $1-\lambda_2>0$ now proves \eqref{e90}, with
$0<\sigma_*\le\frac{1}{2}$ after decreasing it if necessary.
\end{proof}

For fixed $\eps,u$, define
\[
 v_{\eps,u}^0=\frac{\psi_{\eps,u}}{2}(au^{a-1}+bu^{b-1}),\quad
 v_{\eps,u}^1=\frac{\psi_{\eps,u}}{2\eps}(bu^{b-1}-au^{a-1}).
\]
\begin{proposition}\label{p7}
There are $\eps_0>0$ and $0<\sigma\le\frac{1}{2}$, depending only on
$n,\alpha$, such that, for $0<\eps<\eps_0$, $u\in\mathcal{E}_\eps$,
and $z\perp v_{\eps,u}^0,v_{\eps,u}^1,\psi_{\eps,u}\chi_i$,
\[
 \|z\|_2^2-\HLS(\psi_{\eps,u}z,\psi_{\eps,u}z)
 \ge\sigma\|z\|_2^2.
\]
We may take $\eps_0<\min\{q-1,2-q\}$.
\end{proposition}
\begin{proof}
Otherwise take $\eps_j\rightarrow0^+$ as $j\rightarrow\infty$, $u_j\in\mathcal{E}_{\eps_j}$ and
unit vectors $z_j$ satisfying the constraints with deficit
less than $\frac{\sigma_*}{2}$. After taking a subsequence, let $z_j\rightharpoonup z$ as $j\rightarrow\infty$.
Lemma~\ref{l14} implies the limiting orthogonality conditions, while
Lemma~\ref{l13} yields norm convergence to the compact operator
$\mathcal{K}_{\psi_*}$. Hence
\[
 \HLS(\psi_jz_j,\psi_jz_j)\rightarrow\HLS(\psi_*z,\psi_*z)
 \le(1-\sigma_*)\|z\|_2^2\le1-\sigma_*\ \text{as }j\rightarrow\infty,
\]
a contradiction. Shrink $\eps_0$ and take $\sigma=\min\{\frac{\sigma_*}{2},\frac{1}{2}\}$.
\end{proof}

For $1<r<2$, set $B_r(t,s)=t^r-s^r-rs^{r-1}(t-s)$.
Integral Taylor expansion and H\"older yield, for $t\ge0$, $s>0$ and
normalized $f\ge0$, $u\in\mathcal{E}_\eps$,
\begin{align}
 B_r(t,s)&\ge\frac{r(r-1)}{2}\frac{(t-s)^2}{(t+s)^{2-r}},
 \label{e95}\\
 \int_{\R^n} B_r(f,u) dx&\ge
 \frac{r(r-1)}{2(1+2^{\frac{1}{r}})^{2-r}}\|f-u\|_r^2,
 \quad r=a,b.\label{e96}
\end{align}
For \eqref{e96}, use
$\|f-u\|_r^2\le\|f+u\|_r^{2-r}\int_{\R^n} \frac{|f-u|^2}{(f+u)^{2-r}} dx$
and $\|f+u\|_r\le1+2^{\frac{1}{r}}$.

\begin{lemma}\label{l16}
Fix $0<\eps<\eps_0$ and $u\in\mathcal{E}_\eps$. There are
$\eta_u,c_u>0$ such that every normalized $f\ge0$ with
$\Xnorm{f-u}<\eta_u$ and
\begin{equation}\label{e97}
 \int_{\R^n}(f-u)\chi_i dx=0,\quad1\le i\le n,
\end{equation}
satisfies the following remainder estimate
\begin{equation}\label{euclidean-bregman-coercivity}
 C_\eps-\HLS(f,f)\ge c_u\tau_u(f)^2,\quad
 \tau_u(f)^2=\frac{C_\eps}{q}\int_{\R^n}[B_a(f,u)+B_b(f,u)] dx.
\end{equation}
In particular, after adjusting $c_u$,
\begin{equation}\label{e98}
 C_\eps-\HLS(f,f)\ge c_u\Xnorm{f-u}^2.
\end{equation}
\end{lemma}
\begin{proof}
If \eqref{euclidean-bregman-coercivity} fails locally, choose centered
normalized $f_j=u+h_j\rightarrow u$ in $X_\eps$ as $j\rightarrow\infty$, with $h_j\ne0$, and set
$R_j=B_a(f_j,u)+B_b(f_j,u)$, $\tau_j^2=(\frac{C_\eps}{q})\int_{\R^n} R_j dx$.
Then $\tau_j\rightarrow0$ as $j\rightarrow\infty$, $\Xnorm{h_j}\le C\tau_j$ by \eqref{e96}, and
\begin{equation}\label{e99}
 \frac{\delta_j}{\tau_j^2}\rightarrow0,\quad \delta_j=C_\eps-\HLS(f_j,f_j)\ \text{as }j\rightarrow\infty.
\end{equation}
Normalization and \eqref{e35} yield $2\HLS(u,h_j)=-\tau_j^2$, hence
\begin{equation}\label{e100}
 \delta_j=\tau_j^2-\HLS(h_j,h_j).
\end{equation}
To apply Proposition~\ref{p7}, set
$w_j=C_\eps\frac{R_j}{qh_j^2}$,
$\widetilde\psi_j=w_j^{-\frac{1}{2}}$ and
$z_j=\frac{w_j^{\frac{1}{2}}h_j}{\tau_j}$ on $\{h_j\ne0\}$.
Elsewhere set $w_j=W_{\eps,u}$, $\widetilde\psi_j=\psi_{\eps,u}$
and $z_j=0$.
Then
\begin{equation}\label{e101}
 \|z_j\|_2=1,\quad \frac{h_j}{\tau_j}=\widetilde\psi_jz_j.
\end{equation}
After extraction $h_j\rightarrow0$ almost everywhere as $j\rightarrow\infty$. Taylor expansion and
\eqref{e95} yield $w_j\rightarrow W_{\eps,u}$ as $j\rightarrow\infty$ and
$|\widetilde\psi_j|^{s_\alpha}\le C(f_j+u)^q$.
Since $(f_j+u)^q\rightarrow(2u)^q$ in $L^1$ as $j\rightarrow\infty$, Vitali yields
\begin{equation}\label{e102}
 \widetilde\psi_j\rightarrow\psi_{\eps,u}\quad\text{in }L^{s_\alpha}\ \text{as }j\rightarrow\infty.
\end{equation}
For $r=a,b$, \eqref{e95} also yields
$|\widetilde\psi_j u^{r-1}|^2\le C(f_j+u)^{2-r}u^{2r-2}$.
The right side converges in $L^1$ by H\"older and strong $L^r$
convergence: use the conjugate exponents $\frac{r}{2-r}$ and
$\frac{r}{2r-2}$ on $(f_j+u)^{2-r}$ and $u^{2r-2}$.
Pointwise convergence and this $L^1$ control prove
\begin{equation}\label{e103}
 \widetilde\psi_j u^{r-1}\rightarrow\psi_{\eps,u}u^{r-1}
 \quad\text{in }L^2\ \text{as }j\rightarrow\infty.
\end{equation}
Likewise $\widetilde\psi_j\chi_i\rightarrow\psi_{\eps,u}\chi_i$ in $L^2$ as $j\rightarrow\infty$.
These multiplier conclusions use only $f_j\rightarrow u$ as $j\rightarrow\infty$ and the convex
remainders, not \eqref{e99}.

The Euclidean version of \eqref{e40} implies
$|\int_{\R^n} f_j^a dx-1|+|\int_{\R^n} f_j^b dx-1|=O(\sqrt{\delta_j})=o(\tau_j)$.
Since $\int_{\R^n} B_r(f_j,u) dx=O(\tau_j^2)$,
\[
 r\int_{\R^n} u^{r-1}\frac{h_j}{\tau_j} dx
 =\frac{\int_{\R^n} f_j^r dx-1-\int_{\R^n} B_r(f_j,u) dx}{\tau_j}\rightarrow0,
 \quad r=a,b\ \text{as }j\rightarrow\infty.
\]
If $z_j\rightharpoonup z$ as $j\rightarrow\infty$, \eqref{e97}, \eqref{e101} and
\eqref{e103} yield $z\perp a\psi_{\eps,u}u^{a-1}$,
$b\psi_{\eps,u}u^{b-1}$ and $\psi_{\eps,u}\chi_i$.
Since $\eps>0$ is fixed, their half sum and rescaled difference
yield both constraints $z\perp v_{\eps,u}^0,v_{\eps,u}^1$ in
Proposition~\ref{p7}.
Lemma~\ref{l13} and \eqref{e102} therefore yield
\[
 \HLS\Big(\frac{h_j}{\tau_j},\frac{h_j}{\tau_j}\Big)\rightarrow
 \HLS(\psi_{\eps,u}z,\psi_{\eps,u}z)
 \le(1-\sigma)\|z\|_2^2\le1-\sigma\ \text{as }j\rightarrow\infty.
\]
This contradicts \eqref{e99}-\eqref{e100}.
Finally \eqref{e96} implies \eqref{e98}.
\end{proof}

\begin{lemma}\label{l17}
After reducing $\eps_0$, every $u\in\mathcal{E}_\eps$ and every $f$
sufficiently close to $u$ in $L^q$ admit $y(f)\rightarrow0$ as $f\rightarrow u$ such
that $\int_{\R^n}(f(\cdot+y(f))-u)\chi_i dx=0$ for all $i$.
If $f_j\rightarrow u$ in $X_\eps$ as $j\rightarrow\infty$, then $f_j(\cdot+y(f_j))\rightarrow u$ in $X_\eps$ as $j\rightarrow\infty$.
\end{lemma}
\begin{proof}
The map $G_i(f,y)=\int_{\R^n} f(x+y)\chi_i(x) dx$ is $C^1$ on
$L^q\times\R^n$, and $G(u,0)=0$ by radiality. At $(U_*,0)$,
\[
 \partial_{y_j}G_i=-\int_{\R^n} U_*\partial_j\chi_i dx
 =\int_{\R^n} U_*^{q-2}\partial_iU_*\partial_jU_* dx
 =m_*\delta_{ij},\quad m_*>0.
\]
By \eqref{e61}, this matrix remains invertible for all
$u\in\mathcal{E}_\eps$ when $\eps$ is small. The implicit function
theorem provides $y(f)$. The last assertion follows from
$\Xnorm{f_j(\cdot+y_j)-u}\le\Xnorm{f_j-u}+\Xnorm{u(\cdot+y_j)-u}$
and strong continuity of translations.
\end{proof}

\begin{theorem}\label{t10}
For $0<\eps<\eps_0$, there is $k_{\eps,{\rm E}}>0$ such that
\eqref{e18} holds for all normalized nonnegative $f$.
\end{theorem}
\begin{proof}
A sequence contradicting \eqref{e18} has deficit tending to zero,
because $d_{X_{\eps},\rm E}\le2^{\frac{1}{a}}+1<3$.
Theorem~\ref{t6}, the equality classification and Lemma~\ref{l17}
yield translations after which $f_j\rightarrow u\in\mathcal{E}_\eps$ in $X_\eps$ as $j\rightarrow\infty$
and \eqref{e97} holds. Lemma~\ref{l16} yields
$C_\eps-\HLS(f_j,f_j)\ge c_u\Xnorm{f_j-u}^2
\ge c_u(d_{X_{\eps},\rm E}(f_j))^2$, a contradiction.
\end{proof}

\begin{corollary}\label{c2}
The set $\mathcal{E}_\eps$ is finite for $0<\eps<\eps_0$.
\end{corollary}
\begin{proof}
Otherwise radial compactness yields distinct $u_j\rightarrow u$ in
$\mathcal{E}_\eps$. Radiality ensures \eqref{e97}, so
Lemma~\ref{l16} contradicts their equal energies.
\end{proof}

\section{Affine stability}\label{s6}

Christ's shape estimate supplies a near-identity ellipsoid normalization.
We show that the affine correction is $o(\tau^2)$ along any putative
counterexample, where $\tau$ is the scale determined by the convex remainders in Lemma~\ref{l16}.
The Euclidean remainder estimate then yields a contradiction.
Reduce $\eps_0$ so that, throughout this section,
\begin{equation}\label{aff-range}
 1<a=q-\eps<q<b=q+\eps<\min\{2,\kappa\},\quad \kappa=\frac{n}{\alpha}.
\end{equation}
Constants may depend on the fixed positive $\eps$.

We first establish the compactness of the directional bilinear
operators needed in the affine expansion.

Let $S=n\omega_n$, $d\mu=S^{-1}d\xi$, and define the symmetric form
\[
 \mathcal{B}_\alpha(v,h)(\xi)=\frac{1}{2}\int_{\R^n}\int_\R
 |t|^{\alpha-1}v(x+t\xi)h(x)dtdx.
\]
For real $v,h\in L^q$ it is absolutely convergent almost everywhere.
We use the same formula for other inputs whenever the estimates below
ensure convergence. With $D_h=\mathcal{B}_\alpha(h,h)$ for signed $h$,
\begin{equation*}
 D_{u+h}=D_u+2\mathcal{B}_\alpha(u,h)+D_h.
\end{equation*}
Moreover $\HLS(h,h)=S\int_{S^{n-1}} D_hd\mu$ and, for $f\ge0$,
$\A(f)=S\|D_f\|_{L^\kappa(\mu)}$.
By Lemma~\ref{l2}, positivity, and polarization,
\begin{equation}\label{aff-directional-bounds}
 \|D_h\|_\kappa\le C\|h\|_q^2,\quad
 \|\mathcal{B}_\alpha(v,h)\|_\kappa\le C\|v\|_q\|h\|_q.
\end{equation}
Indeed $|D_h|\le D_{|h|}$. Apply the nonnegative bound to
$\lambda|v|+\lambda^{-1}|h|$ and optimize in $\lambda>0$.
Unmarked sphere norms below use $\mu$.

\begin{lemma}\label{aff-ball}
If $1<p\le s<\infty$ and
\begin{equation}\label{aff-cap-condition}
 E(p,s):=\alpha-\frac{1}{p}-\frac{n-1}{s}<0,
\end{equation}
then $\mathcal{B}_\alpha(v,\cdot):L^p\rightarrow L^s(S^{n-1},\mu)$ is
compact for every bounded radial $v$ supported in a ball.
If $|v|\le M1_{B_R}$, its norm is at most $CMR^{n+\alpha-\frac{n}{p}}$.
\end{lemma}
\begin{proof}
For $v=1_{B_1}$, its kernel
$G(x,\xi)=\frac{1}{2}\int_\R|t|^{\alpha-1}1_{B_1}(x+t\xi)dt$
is uniformly bounded on bounded $x$ sets. For $|x|\ge2$, the line
can meet $B_1$ only in two caps of radius $\frac{C}{|x|}$, and $|t|\asymp|x|$
on the segment of intersection. Thus
\begin{equation*}
 G(x,\xi)\le C|x|^{\alpha-1}
 1_{\{\operatorname{dist}(\frac{x}{|x|},\{\xi,-\xi\})\le \frac{C}{|x|}\}}.
\end{equation*}
For $\Omega_R=\{R<|x|<2R\}$, $R\ge2$, radial H\"older and the
cap convolution bound yield the following estimate. Indeed, set
$H_R(\omega)=(\int_R^{2R}|h(r\omega)|^p r^{n-1}dr)^{\frac{1}{p}}$.
The radial integral is bounded by
$CR^{\alpha-1+\frac{n}{p'}}H_R(\omega)$, and
$\|H_R\|_{L^p(S^{n-1},d\omega)}=\|h1_{\Omega_R}\|_p$. Hence
\begin{equation}\label{aff-shell-bound}
 \|\mathcal{B}_\alpha(1_{B_1},h1_{\Omega_R})\|_s
 \le CR^{\alpha-1+\frac{n}{p'}}
 R^{-(n-1)(1+\frac{1}{s}-\frac{1}{p})}\|h1_{\Omega_R}\|_p
 =CR^{E(p,s)}\|h1_{\Omega_R}\|_p.
\end{equation}
The cap operator has $L^p\rightarrow L^p$ norm bounded by its measure and
$L^p\rightarrow L^\infty$ norm bounded by its $\frac{1}{p'}$ power. Interpolation
proves this exponent. The shell estimate holds without a
sign restriction on $E$. Under \eqref{aff-cap-condition}, its dyadic
sum makes the far-field operator norm tend to zero.

On $B_L$, rotational covariance
$G(x,O\xi)=G(O^{-1}x,\xi)$ and strong continuity of rotations show
that $\xi\mapsto G(\cdot,\xi)$ is continuous in $L^{p'}(B_L)$.
The truncated operator therefore maps the $L^p$ unit ball to an
equicontinuous bounded family, compact in $C(S^{n-1})$.
For general radial $v$, the same covariance and domination by a ball
yield the identical argument. Scaling yields the norm bound.
\end{proof}

\begin{lemma}\label{aff-weighted-ball}
Suppose $0\le\beta<n$, $\beta>2\alpha-n$, and
$|\psi(x)|\le C|x|^{-\frac{\beta}{2}}$. For bounded radial $v$ supported in
a ball, $z\mapsto\mathcal{B}_\alpha(v,\psi z)$ is compact $L^2\rightarrow L^2$.
If $|v|\le M1_{B_R}$, its norm is at most
$C'MR^{\alpha+\frac{n}{2}-\frac{\beta}{2}}$.
\end{lemma}
\begin{proof}
For the unit ball, the bounded-$x$ kernel is Hilbert-Schmidt since
$\int_{B_L}|x|^{-\beta} dx<\infty$. Equation \eqref{aff-shell-bound}
with $p=s=2$ bounds the shell operator by
$CR^{\alpha-\frac{n}{2}-\frac{\beta}{2}}$, a negative power of $R$.
Its dyadic tail tends to zero in norm. Domination and scaling yield
the remaining assertions.
\end{proof}

\begin{proposition}\label{aff-directional-compactness}
For $u\in\mathcal{E}_\eps$, there is $s>\max\{\kappa,b\}$ such that
$\mathcal{B}_\alpha(u,\cdot):X_\eps\rightarrow L^s$ is compact. If $\kappa<2$,
one can take $s<2$. Also
$T_uz=\mathcal{B}_\alpha(u,\psi_{\eps,u}z)$ is compact $L^2\rightarrow L^2$.
\end{proposition}
\begin{proof}
For $p=a,b$, $E(p,\kappa)=\frac{\alpha}{n}-\frac{1}{p}<0$, so choose $s>\kappa>b$
close enough to $\kappa$, with $s<2$ when $\kappa<2$.
Radial monotonicity yields $u(\rho)\le C_p\rho^{-\frac{n}{p}}$.
Write $u=\sum_ku_k$, where $u_k=u1_{\{2^k<|x|\le2^{k+1}\}}$.
Lemma~\ref{aff-ball} yields
\begin{equation}\label{aff-radial-sum}
 \|\mathcal{B}_\alpha(u_k,\cdot)\|_{L^p\rightarrow L^s}
 \le C2^{k(n+\alpha-\frac{2n}{p})}.
\end{equation}
Use $p=b$ for $k<0$ and $p=a$ for $k\ge0$. Both geometric sums
converge because $a<q<b$. Thus the operator is a norm-convergent
sum of compact operators on $X_\eps$.

For the weighted assertion,
$\psi_{\eps,u}^2\le C_pu^{2-p}\le C'_p|x|^{-\beta_p}$, where
$\beta_p=\frac{2n}{p}-n$. Condition \eqref{aff-range} yields
$0<\beta_p<n$ and $\beta_p>2\alpha-n$.
Lemma~\ref{aff-weighted-ball} produces the same summable bound
\eqref{aff-radial-sum}, since
$\alpha+\frac{n}{2}-\frac{\beta_p}{2}-\frac{n}{p}=n+\alpha-\frac{2n}{p}$.
These sums agree with the integral formulas by monotone convergence
for nonnegative inputs and absolute domination for signed inputs.
\end{proof}

We next expand the affine correction for perturbations in $X_\eps$,
using the convex remainders to determine their scale.

\begin{proposition}\label{aff-expansion}
Let $u\in\mathcal{E}_\eps$ and let $f_j=u+h_j\ge0$ be nonzero,
with $0\ne h_j\rightarrow0$ in $X_\eps$ as $j\rightarrow\infty$. Set
\begin{gather*}
 R_j=B_a(f_j,u)+B_b(f_j,u),\quad
 \tau_j^2=\frac{C_\eps}{q}\int_{\R^n} R_j dx,
 \\
 w_j=\frac{C_\eps R_j}{qh_j^2},\quad
 z_j=\frac{w_j^{\frac{1}{2}}h_j}{\tau_j}\quad\text{on }\{h_j\ne0\},
\end{gather*}
and $w_j=W_{\eps,u}$, $z_j=0$ elsewhere. If, along a subsequence,
$h_j\rightarrow0$ almost everywhere as $j\rightarrow\infty$ and $z_j\rightharpoonup z$ in $L^2$ as $j\rightarrow\infty$,
then $\|z_j\|_2=1$, $\tau_j\rightarrow0$ as $j\rightarrow\infty$, and, for $d=D_u=\frac{C_\eps}{S}$,
\begin{align}
 \frac{\A(f_j)-\HLS(f_j,f_j)}{\tau_j^2}
 &\rightarrow\frac{2S(\kappa-1)}{d}
 \bigg[\int_{S^{n-1}}|T_uz|^2d\mu-\bigg(\int_{S^{n-1}} T_uzd\mu\bigg)^2\bigg],
 \label{aff-second-order}\\
 \frac{\rho_{K_{f_j}}^n-1}{\tau_j}
 &\rightarrow\frac{2\kappa}{d}\bigg(T_uz-\int_{S^{n-1}} T_uzd\mu\bigg)
 \quad\text{in }L^1(\mu)\ \text{as }j\rightarrow\infty.
 \label{aff-first-order}
\end{align}
\end{proposition}
\begin{proof}
Since $\|f_j+u\|_p$ is bounded, \eqref{e95} and H\"older yield
$\|h_j\|_p^2\le C\int_{\R^n} B_p(f_j,u) dx$ for $p=a,b$.
Thus $\|h_j\|_q=O(\tau_j)$, while continuity of the power integrals
and their linear terms yields $\tau_j\rightarrow0$ as $j\rightarrow\infty$. By
\eqref{aff-directional-bounds},
\begin{equation}\label{aff-directional-orders}
 \|D_{h_j}\|_\kappa=O(\tau_j^2),\quad
 \|\mathcal{B}_\alpha(u,h_j)\|_\kappa=O(\tau_j)\ \text{as }j\rightarrow\infty.
\end{equation}

Fix $0<\eta<\frac{1}{4}$ and split $h_j=g_j+k_j$ with
$g_j=h_j1_{\{|h_j|\le\eta u\}}$.
For $p=a,b$ and $h\ge-u$, $|h|>\eta u$, strict convexity yields
$B_p(u+h,u)\ge c_{p,\eta}|h|^p$: divide by $u^p$ and minimize the
positive ratio on $[-1,-\eta]\cup[\eta,\infty)$, whose limit at
infinity is $1$. Therefore
\begin{equation*}
 \|k_j\|_p\le C\tau_j^{\frac{2}{p}}=o(\tau_j),\quad p=a,b\ \text{as }j\rightarrow\infty.
\end{equation*}
On the other set, Taylor's formula yields
$c_\eta W_{\eps,u}\le w_j\le C_\eta W_{\eps,u}$.
Hence $\frac{g_j}{\tau_j}=\psi_{\eps,u}m_jz_j$, where
$m_j=\Big(\frac{w_j^{-1/2}}{\psi_{\eps,u}}\Big)1_{\{|h_j|\le\eta u\}}$
is uniformly bounded and tends to $1$ almost everywhere, since
$u>0$ almost everywhere. For every $\phi\in L^2$,
$\|(m_j-1)\phi\|_2\rightarrow0$ as $j\rightarrow\infty$, so boundedness of $z_j$ and its weak
convergence yield $m_jz_j\rightharpoonup z$ as $j\rightarrow\infty$.
Proposition~\ref{aff-directional-compactness}
then yields
\begin{align}
 \frac{\mathcal{B}_\alpha(u,g_j)}{\tau_j}&\rightarrow T_uz\quad\text{in }L^2(\mu),
 \label{aff-small-part-limit}\\
 \|\mathcal{B}_\alpha(u,k_j)\|_s
 &\le C(\tau_j^{\frac{2}{a}}+\tau_j^{\frac{2}{b}})=O(\tau_j^{\frac{2}{b}})=o(\tau_j)\ \text{as }j\rightarrow\infty.
 \label{aff-large-part-direction}
\end{align}
Also $|\mathcal{B}_\alpha(u,g_j)|\le\eta d$ by positivity.

Set $e_j=2\mathcal{B}_\alpha(u,g_j)$, $v_j=2\mathcal{B}_\alpha(u,k_j)$,
$\ell_j=e_j+v_j$, and $F(t)=|t|^\kappa$.
We first prove the integrated expansion
\begin{equation}\label{aff-power-expansion}
 \int_{S^{n-1}} F(d+\ell_j)d\mu
 =d^\kappa+\kappa d^{\kappa-1}\int_{S^{n-1}}\ell_j d\mu
 +\frac{\kappa(\kappa-1)}{2}d^{\kappa-2}\int_{S^{n-1}} e_j^2d\mu
 +o(\tau_j^2)\ \text{as }j\rightarrow\infty.
\end{equation}
When $\kappa\ge2$, $s>\kappa$ and
\eqref{aff-large-part-direction} yield $\|v_j\|_2=o(\tau_j)$.
Thus $\frac{\ell_j}{\tau_j}\rightarrow2T_uz$ in $L^2$ as $j\rightarrow\infty$, while
$\|\ell_j\|_\kappa=O(\tau_j)$ by \eqref{aff-directional-orders}.
For $2<\kappa<3$, scalar Taylor expansion at $d$ has integrated
remainder bounded by $C\|\ell_j\|_\kappa^\kappa=o(\tau_j^2)$.
For $\kappa\ge3$, the bound is
$C(\|\ell_j\|_\kappa^\kappa+\|\ell_j\|_\kappa^3)=o(\tau_j^2)$.
For $\kappa=2$, $F$ is quadratic.
In every case
$\int_{S^{n-1}}(\ell_j^2-e_j^2)d\mu=o(\tau_j^2)$ by Cauchy-Schwarz,
which proves \eqref{aff-power-expansion} in this range.

For $1<\kappa<2$, choose $\kappa<s<2$ as above. Since $s>b$,
$\int_{S^{n-1}}|v_j|^s d\mu=O(\tau_j^{2s/b})=o(\tau_j^2)$.
As $|e_j|\le2\eta d<\frac{d}{2}$, scalar Taylor estimates yield
\[
 |F(d+e_j+v_j)-F(d+e_j)-F'(d+e_j)v_j|
 \le C|v_j|^s.
\]
For small $|v_j|$ this uses the second order remainder away from zero.
Otherwise the $\kappa$-power growth suffices.
On $|v_j|\le1$, the mixed term obeys
\[
 |(F'(d+e_j)-F'(d))v_j|
 \le\delta|e_j|^2+C_\delta|v_j|^s\quad(\delta>0),
\]
and on $|v_j|>1$ it is bounded by $C|v_j|^s$.
For fixed $\delta>0$, integration therefore yields
\[
 \limsup_{j\rightarrow\infty}\frac{1}{\tau_j^2}
 \int_{S^{n-1}}|(F'(d+e_j)-F'(d))v_j|d\mu
 \le\delta\limsup_{j\rightarrow\infty}\Big\|\frac{e_j}{\tau_j}\Big\|_2^2.
\]
As $\delta\rightarrow0^+$, we conclude that this mixed term is
$o(\tau_j^2)$. For the term involving only $e_j$, write
\[
 F(d+e_j)-F(d)-F'(d)e_j-\frac{1}{2}F''(d)e_j^2
 =e_j^2\omega_d(e_j),
\]
where $\omega_d$ is bounded on $[-\frac{d}{2},\frac{d}{2}]$ and
$\omega_d(t)\rightarrow0$ as $t\rightarrow0$. By
\eqref{aff-small-part-limit}, $\frac{e_j}{\tau_j}$ converges strongly in
$L^2$, so $e_j^2/\tau_j^2$ is uniformly integrable and $e_j\rightarrow0$ as $j\rightarrow\infty$
in measure. Splitting into $|e_j|\le r$ and $|e_j|>r$, first
as $j\rightarrow\infty$ and then as $r\rightarrow0^+$, proves
$\int_{S^{n-1}} e_j^2|\omega_d(e_j)|d\mu=o(\tau_j^2)$.
This proves \eqref{aff-power-expansion} also when $1<\kappa<2$.
No $L^2$ bound on $v_j$ has been used in this range.

Let $N(w)=\|w\|_\kappa$, $M_j=\int_{S^{n-1}}\ell_j d\mu$ and
$Q_j=\int_{S^{n-1}} e_j^2d\mu$. Since $M_j=O(\tau_j)$ and
$Q_j=O(\tau_j^2)$, scalar expansion of the outer power in
\eqref{aff-power-expansion} yields
\[
 N(d+\ell_j)=d+M_j+\frac{\kappa-1}{2d}(Q_j-M_j^2)
 +o(\tau_j^2)\ \text{as }j\rightarrow\infty.
\]
Also $N$ is continuously Fr\'echet differentiable on
$L^\kappa\setminus\{0\}$, with
$N'(w)=\frac{|w|^{\kappa-2}w}{\|w\|_\kappa^{\kappa-1}}$ as a dual
density and $N'(d)=1$. Set $r_j=D_{h_j}$.
Since $\|\ell_j\|_\kappa=O(\tau_j)$ and
$\|r_j\|_\kappa=O(\tau_j^2)$, the fundamental theorem of calculus
yields
\[
 \bigg|N(d+\ell_j+r_j)-N(d+\ell_j)-\int_{S^{n-1}} r_j d\mu\bigg|
 \le\|r_j\|_\kappa
 \sup_{0\le t\le1}\|N'(d+\ell_j+tr_j)-1\|_{\kappa'}
 =o(\tau_j^2)\ \text{as }j\rightarrow\infty.
\]
All arguments of $N'$ stay away from zero in norm for large $j$.
Finally $\int_{S^{n-1}} v_j d\mu=o(\tau_j)$ implies
$M_j^2=(\int_{S^{n-1}} e_j d\mu)^2+o(\tau_j^2)$. Subtracting
$\int_{S^{n-1}} D_{f_j}d\mu=d+M_j+\int_{S^{n-1}} r_j d\mu$ yields
\begin{equation*}
 \|D_{f_j}\|_\kappa-\int_{S^{n-1}} D_{f_j}d\mu
 =\frac{\kappa-1}{2d}
 \bigg[\int_{S^{n-1}} e_j^2d\mu-\bigg(\int_{S^{n-1}} e_j d\mu\bigg)^2\bigg]
 +o(\tau_j^2)\ \text{as }j\rightarrow\infty.
\end{equation*}
This proves \eqref{aff-second-order}.

For \eqref{aff-first-order}, \eqref{e63} yields
\begin{equation}\label{aff-star-density}
 \rho_{K_f}^n=\frac{D_f^\kappa}{\int_{S^{n-1}} D_f^\kappa d\mu}.
\end{equation}
The map $v\mapsto\frac{|v|^\kappa}{\int_{S^{n-1}}|v|^\kappa d\mu}$ is Fr\'echet
differentiable $L^\kappa\rightarrow L^1$ at $d$, with derivative
$k\mapsto(\frac{\kappa}{d})(k-\int_{S^{n-1}} kd\mu)$.
Indeed the power map has derivative
$k\mapsto\kappa d^{\kappa-1}k$ from $L^\kappa$ to $L^1$.
Its remainder is bounded in $L^1$ by $C\|k\|_\kappa^\kappa$
for $1<\kappa\le2$, and by
$C(\|k\|_\kappa^2+\|k\|_\kappa^\kappa)$ for $\kappa>2$.
Dividing by its integral yields the stated normalized derivative.
Here $k=D_{f_j}-d=\ell_j+r_j$ has norm $O(\tau_j)$, so the
derivative remainder is $o(\tau_j)$.
Equations \eqref{aff-small-part-limit}-\eqref{aff-large-part-direction}
and \eqref{aff-directional-orders} yield $\frac{k}{\tau_j}\rightarrow2T_uz$ in
$L^1$, proving \eqref{aff-first-order}.
\end{proof}

The preceding expansions will be combined with a linear change
of variables that controls the associated star set. We first record
the geometric identities needed for this normalization.

Polar coordinates and changes of variables yield
\begin{equation}\label{aff-shape-identity}
 \|\rho_K^n-\rho_E^n\|_{L^1(\mu)}
 =\frac{|K\mathbin\triangle E|}{\omega_n},\quad
 K_{f\circ L^{-1}}=LK_f\quad(L\in\SL(n)).
\end{equation}
For the covariance identity, explicitly,
\[
 D_{f\circ L^{-1}}(\xi)=|L^{-1}\xi|^{-\alpha}
 D_f\bigg(\frac{L^{-1}\xi}{|L^{-1}\xi|}\bigg).
\]
This yields $S_\alpha(f\circ L^{-1})=LS_\alpha f$.
Volume normalization commutes with $L$ because $\det L=1$.
Translations leave $D_f$ and $K_f$ unchanged.

\begin{lemma}\label{aff-ellipsoid-normalization}
If normalized $f_j\ge0$ converge to $u\in\mathcal{E}_\eps$ in $X_\eps$
and $\Delta_j=\Delta_\eps(f_j)>0$ tends to zero, there are
$L_j\in\SL(n)$ with $L_j\rightarrow I$ as $j\rightarrow\infty$ such that $g_j=f_j\circ L_j^{-1}$
satisfies
\begin{equation*}
 g_j\rightarrow u\text{ in }X_\eps,\quad
 \|\rho_{K_{g_j}}^n-1\|_1\le C\sqrt{\Delta_j}\ \text{as }j\rightarrow\infty.
\end{equation*}
\end{lemma}
\begin{proof}
By \eqref{aff-directional-bounds} and \eqref{aff-star-density},
$\rho_{K_{f_j}}^n\rightarrow1$ in $L^1$ as $j\rightarrow\infty$.
Theorem~\ref{t7} yields $\mathfrak{a}(K_{f_j})\le C\sqrt{\Delta_j}$.
Choose centered volume-normalized ellipsoids $E_j=A_jB_1$, with
$A_j$ symmetric positive definite and $\det A_j=1$, such that
\[
 |K_{f_j}\mathbin\triangle E_j|
 \le\omega_n\mathfrak{a}(K_{f_j})+\sqrt{\Delta_j}.
\]
Then $E_j\rightarrow B_1$ in measure as $j\rightarrow\infty$. The least eigenvalue of $A_j$ cannot
tend to zero, since $E_j\cap B_1$ would lie in a slab of vanishing
thickness. The determinant condition bounds the largest eigenvalue.
Every subsequential matrix limit is positive definite and maps $B_1$
onto itself, hence is $I$. Thus $A_j\rightarrow I$ as $j\rightarrow\infty$.
Take $L_j=A_j^{-1}$. Then \eqref{aff-shape-identity} yields
\[
 \|\rho_{K_{g_j}}^n-1\|_1
 =\frac{|L_jK_{f_j}\triangle B_1|}{\omega_n}
 =\frac{|K_{f_j}\triangle E_j|}{\omega_n}
 \le C\sqrt{\Delta_j}.
\]
Finally $\Xnorm{f_j\circ A_j-u}\le\Xnorm{f_j-u}+\Xnorm{u\circ A_j-u}\rightarrow0$ as $j\rightarrow\infty$.
Strong continuity of the action follows by approximation with smooth
compactly supported functions.
\end{proof}

We combine the ellipsoid normalization, the nonlinear affine
expansion and the Euclidean estimate to prove Theorem~\ref{t2}.
\begin{proof}
Theorem~\ref{t10} and Corollary~\ref{c2} yield the Euclidean estimate
and finiteness. If the affine estimate fails, choose normalized $f_j\ge0$
with $d_j=d_{X_{\eps},{\rm A}}(f_j)>0$ and $\frac{\Delta_j}{d_j^2}\rightarrow0$ as $j\rightarrow\infty$.
Since $d_j\le2^{\frac{1}{a}}+1<3$, their deficits tend to zero and are positive
by Theorem~\ref{t8}. Theorem~\ref{t9} supplies affine transformations
and translations after which $f_j\rightarrow u\in\mathcal{E}_\eps$ in $X_\eps$ as $j\rightarrow\infty$.
Apply Lemma~\ref{aff-ellipsoid-normalization}, then the translation
centering of Lemma~\ref{l17}. The resulting sequence still converges
to $u$, satisfies \eqref{e97}, and obeys
$\|\rho_{K_{f_j}}^n-1\|_1\le C\sqrt{\Delta_j}$.
Translations preserve the star set, and all transformations used
are isometries of $X_\eps$ carrying $\mathcal{M}_{\eps,\rm A}$ onto
itself. They therefore preserve $d_j$, $\Delta_j$ and normalization.

Define the variables associated with the convex remainders only after this normalization, as in
Proposition~\ref{aff-expansion}. Equation \eqref{e96} yields
\begin{equation}\label{aff-counter-scale}
 d_j\le\Xnorm{f_j-u}\le C\tau_j,\quad \Delta_j=o(\tau_j^2)\ \text{as }j\rightarrow\infty.
\end{equation}
Along a subsequence, $h_j=f_j-u\rightarrow0$ almost everywhere
and $z_j\rightharpoonup z$ in $L^2$ as $j\rightarrow\infty$.
By \eqref{aff-counter-scale},
$\|\rho_{K_{f_j}}^n-1\|_1\le C\sqrt{\Delta_j}=o(\tau_j)$.
Thus \eqref{aff-first-order} forces $T_uz$ to be constant,
so the variance in \eqref{aff-second-order} vanishes. Hence
\begin{equation*}
 \A(f_j)-\HLS(f_j,f_j)=o(\tau_j^2)\ \text{as }j\rightarrow\infty.
\end{equation*}
Together with \eqref{aff-counter-scale}, this says
$C_\eps-\HLS(f_j,f_j)=o(\tau_j^2)$, contradicting the centered
remainder estimate \eqref{euclidean-bregman-coercivity}.
The $L^q$ estimates follow from $\|h\|_q\le\Xnorm{h}$, and
Proposition~\ref{aff-sharp-power} proves sharpness of the power.
\end{proof}

\section{Sharpness and parameter dependence}\label{s7}

Normalized dilations determine the sharp power and force the
stability constants to vanish as $\eps\rightarrow0^+$.

\begin{proposition}\label{p8}
For $u\in\mathcal{E}_\eps$, set
$s_\eps(t)=\Scale{\mathcal{D}_{e^t}u}$,
$F_{\eps,t}=s_\eps(t)^{-1}\mathcal{D}_{e^t}u$ and
$y_\eps(t)=\log(s_\eps(t))$.
The normalization is independent of $u$ and satisfies
\begin{gather}
qy_\eps(t)=\log\Big(\cosh\Big(\frac{n\eps}{q}t-\eps y_\eps(t)\Big)\Big),
\label{e104}\\
\Delta_\eps(F_{\eps,t})=C_\eps(1-e^{-2y_\eps(t)}).\label{e105}
\end{gather}
Moreover, $s_\eps(0)=1$, $s_\eps(t)>1$ for $t\ne0$, and
\begin{align}
\Delta_\eps(F_{\eps,t})
&=\frac{C_\eps n^2\eps^2}{q^3}t^2+O_\eps(|t|^3)
\ \text{as }t\rightarrow0,\label{e106}\\
\frac{\Delta_\eps(F_{\eps,t})}{C_\eps}
&=\frac{n^2}{q^3}\eps^2t^2+O(\eps^4)
\ \text{as }\eps\rightarrow0^+,\label{e107}
\end{align}
the last expansion being uniform on bounded $t$ intervals.
\end{proposition}
\begin{proof}
Theorem~\ref{t4} and critical scaling yield
\[
1=\frac{1}{2}e^{-n\eps \frac{t}{q}}s_\eps(t)^{-(q-\eps)}
 +\frac{1}{2}e^{n\eps \frac{t}{q}}s_\eps(t)^{-(q+\eps)},
\]
equivalently \eqref{e104}. The function $y\mapsto qy-\log(\cosh(\frac{n\eps}{q}t-\eps y))$ has derivative at least $q-\eps>0$, is nonpositive at zero, and
tends to infinity. This proves uniqueness and the assertions about
$s_\eps$. Energy invariance yields \eqref{e105}.
Implicit differentiation yields $y_\eps'(0)=0$ and
$y_\eps''(0)=\frac{n^2\eps^2}{q^3}$, proving \eqref{e106}.
For bounded $t$,
\[
0\le y_\eps(t)\le
\frac{\log\big(\cosh\big(n\eps \frac{t}{q}\big)\big)}{q-\eps}=O(\eps^2)\ \text{as }\eps\rightarrow0^+.
\]
Substitution into \eqref{e104} yields
$y_\eps(t)=\frac{n^2\eps^2t^2}{2q^3}+O(\eps^4)$, and hence
\eqref{e107}.
\end{proof}

For $p>0$, define the sharp constants
\begin{align}
\mathfrak{c}_{\eps,p}
&=\inf_{\substack{f\ge0,\ \M(f)=1\\d_{q,{\rm A}}(f)>0}}
  \frac{\Delta_\eps(f)}{d_{q,{\rm A}}(f)^p},\label{e108}\\
\mathfrak{c}_{\eps,p,X_{\eps}}
&=\inf_{\substack{f\ge0,\ \M(f)=1\\d_{X_{\eps},{\rm A}}(f)>0}}
  \frac{\Delta_\eps(f)}{d_{X_{\eps},{\rm A}}(f)^p}.\notag
\end{align}
Theorem~\ref{t2} and $d_{q,{\rm A}}\le d_{X_{\eps},{\rm A}}$ imply
\begin{equation*}
0<\mathfrak{c}_{\eps,2,X_{\eps}}\le\mathfrak{c}_{\eps,2}
\quad(0<\eps<\eps_0).
\end{equation*}
Since both distances are at most $2^{\frac{1}{a}}+1<3$, the constants are
also positive for every $p\ge2$.

\begin{proposition}\label{aff-sharp-power}
For $0<\eps<q-1$ and $0<p<2$, $\mathfrak{c}_{\eps,p,X_{\eps}}=0$.
For $0<\eps<\eps_0$, also $\mathfrak{c}_{\eps,p}=0$ for $0<p<2$.
The power $2$ is sharp in all four estimates of Theorem~\ref{t2}.
\end{proposition}
\begin{proof}
Fix $u\in\mathcal{E}_\eps$. Every $v\in\mathcal{M}_{\eps,\rm A}$
has $\|v\|_a=1$. Hence
\begin{align*}
d_{X_{\eps},{\rm A}}(F_{\eps,t})
\ge\big|\|F_{\eps,t}\|_a-1\big|=\Big|e^{-n\eps \frac{t}{qa}-y_\eps(t)}-1\Big|
=\frac{n\eps}{qa}|t|+O_\eps(t^2)\ \text{as }t\rightarrow0.
\end{align*}
Together with \eqref{e106}, this proves the $X_\eps$ assertion.

For small $\eps$, $\mathcal{E}_\eps$ is finite by Corollary~\ref{c2}.
Rearrangement contraction yields, for radial nonincreasing $F$,
\begin{equation}\label{aff-radial-distance}
d_{q,{\rm A}}(F)=\min_{v\in\mathcal{E}_\eps}\|F-v\|_q.
\end{equation}
As $F_{\eps,t}\rightarrow u$ in $L^q$, its nearest member is $u$ for
small $|t|$. The distribution $\Lambda u=(\frac{n}{q})u+x\cdot\nabla u$
is nonzero: otherwise differentiation against compactly supported
tests yields
\[
 \frac{d}{dt}\langle\mathcal{D}_{e^t}u,\phi\rangle
 =\langle\mathcal{D}_{e^t}\Lambda u,\phi\rangle=0.
\]
Thus the dilation orbit would be constant as a distribution and
hence in $L^q$. All disjoint dyadic dilates of
$\{1<|x|<2\}$ would then have the same $L^q$ mass. Integrability
forces each mass to vanish, and their union forces $u=0$.
Choose $\phi\in C_c^\infty$ with
$\langle\Lambda u,\phi\rangle\ne0$.
Differentiating
$\int_{\R^n}\mathcal{D}_{e^t}u\,\phi
=e^{-\frac{nt}{q'}}\int_{\R^n} u(x)\phi(e^{-t}x) dx$
is justified by local integrability. Since $y_\eps'(0)=0$,
\[
\|F_{\eps,t}-u\|_q
\ge\frac{|\langle F_{\eps,t}-u,\phi\rangle|}{\|\phi\|_{q'}}
=\frac{|\langle\Lambda u,\phi\rangle|}{\|\phi\|_{q'}}|t|+o(|t|)\ \text{as }t\rightarrow0.
\]
Equations \eqref{aff-radial-distance} and \eqref{e106} prove the
$L^q$ assertion without assuming strong differentiability of $u$.
Finally, the test functions are radial, so the two deficits agree,
while the Euclidean distances dominate the affine distances.
The same tests exclude $0<p<2$ in the Euclidean estimates.
\end{proof}

\begin{theorem}
For every fixed $p>0$, $\mathfrak{c}_{\eps,p}=O(\eps^2)$ as
$\eps\rightarrow0^+$. More precisely, for $t\ne0$ and
$d_t=\|\mathcal{D}_{e^t}U_*-U_*\|_q>0$,
\begin{equation}\label{e109}
\limsup_{\eps\rightarrow0^+}\frac{\mathfrak{c}_{\eps,p}}{\eps^2}
\le\frac{C_0n^2t^2}{q^3d_t^p}.
\end{equation}
In particular,
\begin{equation}\label{e110}
\limsup_{\eps\rightarrow0^+}\frac{\mathfrak{c}_{\eps,2}}{\eps^2}
\le\frac{C_0}{q}
\bigg(\frac{\sqrt\pi\,\Gamma(\frac{n+q+1}{2})}{\Gamma(\frac{q+1}{2})\Gamma(\frac{n+1}{2})}\bigg)^{\frac{2}{q}}.
\end{equation}
No positive stability constant can be uniform in small $\eps$.
\end{theorem}
\begin{proof}
Choose $u_\eps\in\mathcal{E}_\eps$. For fixed $t\ne0$,
$F_{\eps,t}\rightarrow\mathcal{D}_{e^t}U_*$ in $L^q$ as $\eps\rightarrow0^+$ by Theorem~\ref{t1}
and Proposition~\ref{p8}. Equation \eqref{aff-radial-distance} and
the uniform convergence \eqref{e61} imply
$d_{q,{\rm A}}(F_{\eps,t})\rightarrow d_t>0$ as $\eps\rightarrow0^+$.
Testing \eqref{e108} with $F_{\eps,t}$ and using \eqref{e107}
proves \eqref{e109}.

The explicit dilation orbit of $U_*$ is strongly differentiable,
so $\frac{d_t}{|t|}\rightarrow\|\Lambda U_*\|_q=\|\Lambda U\|_q$.
Since
$\Lambda U=\gamma \frac{U(1-|x|^2)}{1+|x|^2}$, the beta integral yields
\[
\|\Lambda U\|_q^q
=\gamma^q\frac{\Gamma(\frac{q+1}{2})\Gamma(\frac{n+1}{2})}{\sqrt\pi\,\Gamma(\frac{n+q+1}{2})}.
\]
Using $\gamma=\frac{n}{q}$ in \eqref{e109} as $t\rightarrow0$, we obtain
\eqref{e110}.
\end{proof}

\begin{remark}
For fixed small $\eps$,
$0<\mathfrak{c}_{\eps,2,X_{\eps}}\le\mathfrak{c}_{\eps,2}=O(\eps^2)$.
The $X_\eps$ constant satisfies the explicit upper bound
\eqref{e110} as well. The same radial dilation tests yield
$O(\eps^2)$ upper bounds for the sharp Euclidean constants in
both distances. A matching lower bound
$\mathfrak{c}_{\eps,2}\ge c_{n,\alpha}\eps^2$, and the stronger
version with $\mathfrak{c}_{\eps,2,X_{\eps}}$, remain open.
\end{remark}

\section*{Acknowledgments}
This work was supported by
Natural Science Basic Research Program of Shaanxi Province of China (2026JC-YBMS-0065).

\bigbreak

\begin{samepage}
\noindent Zhenzhen Wei\\
Department of Mathematics and Statistics,
Memorial University of Newfoundland,
St. John's, Newfoundland and Labrador A1C 5S7, Canada\\
E-mail: {\tt zhenzhenw@mun.ca}\par\medskip

\noindent Denghui Wu\\
College of Science, Northwest A\&F University,
Yangling, Shaanxi 712100, China\\
E-mail: {\tt wudenghui66@163.com}

\end{samepage}

\end{document}